\documentclass[11pt]{amsart}

\usepackage{amsmath}
\usepackage{amssymb}
\usepackage{graphicx}

\usepackage{dsfont}
\usepackage{enumitem}
\usepackage{tikz}
\usepackage{subcaption}
\usepackage{needspace}

\newtheorem{theorem}{Theorem}[section]
\newtheorem{corollary}[theorem]{Corollary}
\newtheorem{lemma}[theorem]{Lemma}

\theoremstyle{definition}

\newtheorem{remark}[theorem]{Remark}
\newtheorem{assumption}[theorem]{Assumption}

\numberwithin{equation}{section}

\newcommand{\mexp}[2][]{\mathbb E_{#1}\! \left[#2\right]}
\newcommand{\Ai}[1]{\,\mathrm{Ai} \!\left( #1 \right) }
\newcommand{\AiPrime}[1]{\,\mathrm{Ai}' \!\left( #1 \right) }
\newcommand{\Tr}[1]{\mathrm{Tr}\, #1}
\renewcommand{\Re}[1]{\mathrm{Re}\, #1}
\renewcommand{\Im}[1]{\mathrm{Im}\, #1}

\makeatletter
\def\namedlabel#1#2{\begingroup
  \def\@currentlabel{#2}%
  \label{#1}\endgroup}
\makeatother

\newtheoremstyle{prbl}
{}                                   
{}                                   
{\itshape}                           
{}                                   
{\bfseries}                          
{}                                   
{ }                                  
{\thmname{#1}\thmnote{ #3}}          
\theoremstyle{prbl}
\newtheorem{problem}{Problem}

\usetikzlibrary{arrows,decorations.markings}

\title[Rate of Convergence in CLT for Linear Statistics]{On the Rate of Convergence in the Central Limit Theorem for Linear Statistics of Gaussian, Laguerre, and Jacobi Ensembles}

\author[Sergey Berezin]{Sergey Berezin}
\address[Sergey Berezin]{Aix-Marseille Universit\'e, Centrale Marseille, CNRS, Institut de Math\'ematiques de Marseille, UMR7373, 39 Rue F. Joliot Curie 13453, Marseille, France; St. Petersburg Department of V.A.~Steklov Mathematical Institute of RAS, 27 Fontanka 191023, St. Petersburg, Russia}
\email{{\tt servberezin@yandex.ru, sergey.berezin@univ-amu.fr}}

\author[Alexander I. Bufetov]{Alexander I. Bufetov}
\address[Alexander I. Bufetov]{Aix-Marseille Universit\'e, Centrale Marseille, CNRS, Institut de Math\'ematiques de Marseille, UMR7373, 39 Rue F. Joliot Curie 13453, Marseille, France; Steklov Mathematical Institute of RAS, 8 Gubkina 119991, Moscow, Russia; Institute for Information Transmission Problems, 19 Bolshoy Karetny 127994, Moscow, Russia}
\email{\tt bufetov@mi.ras.ru, alexander.bufetov@univ-amu.fr}



\begin{document}
\begin{abstract}
  Under the Kolmogorov--Smirnov metric, an upper bound on the rate of convergence to the Gaussian distribution is obtained for linear statistics of the matrix ensembles in the case of the Gaussian, Laguerre, and Jacobi weights. The main lemma gives an estimate for the characteristic functions of the linear statistics; this estimate is uniform over the growing interval. The proof of the lemma relies on the Riemann--Hilbert approach.
\end{abstract}

\maketitle


\section{Introduction}
\label{subsec_intro}
The main result of this paper is an upper bound, under the Kolmogorov--Smirnov metric, on the rate of convergence in the central limit theorem (CLT) for linear statistics of the Gaussian (GUE), Laguerre (LUE), and Jacobi (JUE) unitary ensembles. Those are canonical examples of the matrix ensembles corresponding to one-cut regular potentials (see~\cite{Charlier2018, Charlier2019}).

We begin by defining the objects of our interest. Let~$\mathcal{M}$ be the set of~$n \times n$ Hermitian matrices, and let~$\mathcal{M}_{\mathcal{I}} \subset \mathcal{M}$ be the set of matrices with eigenvalues in~$\mathcal{I}$, where~$\mathcal{I}$ is the interval given by
\begin{equation}
  \mathcal{I} = \left\{
    \begin{aligned}
      &(-\infty,+\infty)  &&\mbox{for GUE},\\
      &[-1,+\infty) &&\mbox{for LUE},\\
      &[-1,1] &&\mbox{for JUE}.
    \end{aligned}
  \right.
\end{equation}
Endow~$\mathcal{M}$ with the probability measure
\begin{equation}
  \label{nu-measure}
  \mathbb{P}_{n}(dM) = \frac{\mathds{1}_{\mathcal{M}_{\mathcal{I}}}(M)}{Z_{n}} e^{-\Tr{Q_n(M)}} dM,\quad M = \{M_{j,k}\}_{j,k=1}^n \in \mathcal{M},
\end{equation}
where~$dM = \prod\limits_{j} dM_{j,j}  \prod\limits_{j<k}dM^{\mathrm{Re}}_{j,k}\, dM^{{\mathrm{Im}}}_{j,k}$ is the Lebesgue measure on elements of~$M$, and~$Z_{n}$ is the corresponding normalizing constant; the function~$Q_n(x)$ in the exponent is given by
\begin{equation}
  \label{qn-weight}
  Q_n(x) = n V(x) - \omega(x),\quad x \in \mathcal{I};
\end{equation}
the potential~$V(x)$ and additional term~$\omega(x)$ are given respectively by
\begin{equation}
  \label{V_omega}
  V(x) = \left\{
    \begin{aligned}
      &2 x^2  &&\mbox{for GUE},\\
      &2(x +1) &&\mbox{for LUE},\\
      &0 &&\mbox{for JUE},
    \end{aligned}
  \right.
\end{equation}
and
\begin{equation}
  \label{Omega_omega}
  \omega(x) = \left\{
    \begin{aligned}
      &0  &&\mbox{for GUE},\\
      &\alpha \log{(1+x)} &&\mbox{for LUE},\\
      &\alpha \log{(1+x)} + \beta \log{(1-x)} &&\mbox{for JUE},
    \end{aligned}
  \right.
\end{equation}
where~$\alpha, \beta>-1$. Further, we also use~$M$ to denote the random matrix corresponding to~\eqref{nu-measure}, a random element on~$\mathcal{M}$.

The choice of~$\mathcal{I}$, $V(x)$, and~$\omega(x)$ for LUE is somewhat non-standard; nevertheless, it ensures that the corresponding equilibrium measures (see Section~\ref{subsect_eq_measure}) are supported on the same interval~$[-1,1]$ in all three cases (GUE, LUE, and JUE). This convention, also adopted in~\cite{Charlier2019} by Charlier and Gharakhloo, makes it easier to compare the ensembles to each other.

Let~$f$ be a real-valued function on $\mathcal{I}$. Define the linear functionals~$\varkappa[\cdot]$ and $\mu[\cdot]$, respectively, by
\begin{equation}
  \label{eq:ch1_leading_asymp_kappa}
  \varkappa[f] = \left\{
    \begin{aligned}
      &\frac{2}{\pi} \int \limits_{-1}^1 f(x) \sqrt{1-x^2}\, dx  &&\mbox{for GUE},\\
      &\frac{1}{\pi} \int \limits_{-1}^1 f(x) \sqrt{\frac{1-x}{1+x}}\, dx &&\mbox{for LUE},\\
      &\frac{1}{\pi} \int \limits_{-1}^1 f(x) \frac{1}{\sqrt{1-x^2}}\, dx  &&\mbox{for JUE},
    \end{aligned}
  \right.
\end{equation}
and
\begin{equation}
  \mu[f] = \left\{
    \begin{aligned}
      &0  &&\mbox{for GUE},\\
      &\frac{\alpha}{2\pi} \int \limits_{-1}^1 \frac{f(x)-f(-1)}{\sqrt{1-x^2}}\, dx &&\mbox{for LUE},\\
      &\frac{\alpha}{2\pi} \int \limits_{-1}^1 \frac{f(x)-f(-1)}{\sqrt{1-x^2}}\, dx +\frac{\beta}{2\pi} \int \limits_{-1}^1 \frac{f(x)-f(1)}{\sqrt{1-x^2}}\, dx &&\mbox{for JUE};
    \end{aligned}
  \right.  
\end{equation}
also, define the non-negative quadratic functional~$K[\cdot]$ by
\begin{equation}
  \label{variance_kf}
  K[f] = \frac{1}{2 \pi^2} \int \limits_{-1}^1 \frac{f(x)}{\sqrt{1-x^2}} \mathrm{v.p.} \int \limits_{-1}^1\frac{f'(y)\sqrt{1-y^2}}{x-y} \, dy \, dx.
\end{equation}

Let~$F_{f,n}(x) = \mathbb{P}_n\{S_{f,n} \le x\}$ stand for the cumulative distribution function, under the measure \eqref{nu-measure}, of the random variable~$S_{f,n}$,
\begin{equation}
  \label{eq_eq_Sfn}
  S_{f,n}=\frac {\Tr{f(M)} - n \varkappa[f] - \mu[f]}{\sqrt{K[f]}}, 
\end{equation}
and let~$F_{\mathcal N}$ stand for the the cumulative distribution function of the standard Gaussian law of expectation zero and variance one,
\begin{equation}
  F_{\mathcal N} (x) = \frac{1}{\sqrt{2 \pi}} \int \limits_{-\infty}^{x} e^{-s^2/2}\, ds.
\end{equation}

In this setup we have our main result, a theorem that gives an upper bound for the Kolmogorov--Smirnov distance~$\sup \limits_x|F_{f,n}(x)-F_{\mathcal N}(x)|$.
\begin{theorem}
  \label{th1_speed_of_conv}
  Let~$f:\mathcal{I}\to \mathbb R$ be a locally H\"{o}lder continuous function  admitting an analytic continuation into a complex neighborhood of~$[-1,1]$. Additionally for GUE and LUE, let~$f$ satisfy~$f(x) = O(e^{A\, V(x)})$ as~$|x| \to +\infty$, $x \in \mathcal{I}$, for some~$A>0$. Then
  \begin{equation}
    \label{th1_eq1}
    \sup\limits_{n, x}{(n^{1/d}|F_{f,n}(x)-F_{\mathcal N}(x)|)}<+\infty,
  \end{equation}
  where~$d=5$ for GUE and LUE, and~$d=3$ for JUE.
\end{theorem}

\begin{remark}
  The assumption of~$f$ being locally H\"{o}lder continuous function outside the neighborhood of~$[-1,1]$ is technical and can be weakened since the exponential bounds on the tail probabilities for the extreme eigenvalues are known (see~\cite[Proposition~2.1]{Borot2013}). Nevertheless, we do not discuss this in the present manuscript: our choice of the class of test functions is solely governed by the method we use further.  
\end{remark}
\begin{remark}
  For LUE and JUE, the rate of convergence prescribed by this theorem is~$O(n^{-1/d})$, and the latter is uniform with respect to~$\alpha, \beta > -1$ in compact sets. Yet the question remains open to find the rate-optimal estimates and to study their dependence on the smoothness of the test functions. 
\end{remark}
\begin{remark}
  \label{rem_dep_on_loc_param}
  By the very definition of~$S_{f,n}$, its distribution depends only on the distribution of the eigenvalues of~$M$. But it is the local behavior of these eigenvalues near the (asymptotic) edges of the spectrum, the points~$\pm 1$, what determines the rate of convergence in Theorem~\ref{th1_speed_of_conv}. Further discussion will follow in Section~\ref{subsect:outl}.
\end{remark}

A problem similar to the one considered here has also been addressed by several other authors. The rate of convergence in CLT for the traces of powers of the random matrices from compact classical groups was studied by Stein~\cite{Stein1995} and Johansson~\cite{Johansson1997}. In~\cite{Stein1995}, the super-polynomial convergence is proven for the circular real ensemble, which corresponds to the normalized Haar measure on the orthogonal group. Johansson~\cite{Johansson1997} obtained the super-exponential rate of convergence in CLT for the circular unitary ensemble, which corresponds to the normalized Haar measure on the unitary group, and the exponential rate for the circular real and quaternion ensembles, which correspond respectively to the normalized Haar measure on the orthogonal and symplectic groups. Johansson's proof is based on the explicit representation of the moments through combinatorial identities for Toeplitz determinants, on the super-exponential bound for the characteristic function (see \cite[Propositions 2.10, 3.8]{Johansson1997}), and on the classical smoothing inequality of Feller (e.g., see~\cite{Feller2}). In~\cite{Keating2010}, a similar problem to that in~\cite{Johansson1997} was studied for the traces of random matrices from the circular unitary ensemble, prefactored by a deterministic complex matrix. In~\cite{Lambert2019}, Lambert, Ledoux, and Webb considered $\beta$-ensembles with one-cut potentials on the real line. They studied the quadratic Kantorovich distance between the standard Gaussian law and the law of a smooth linear statistic.

In our proof of Theorem~\ref{th1_speed_of_conv}, we also rely on the smoothing inequality of Feller. In order to use this inequality effectively, we need good control over the behavior of the characteristic function of our (asymptotically centered and normalized) linear statistic~$S_{f,n}$. Let~$\varphi_{f,n}(h) = \mexp[n]{e^{i h S_{f,n}}}$, $h \in \mathbb{R}$, be the characteristic function of this statistic, and let ~$\varphi_{\mathcal{N}}(h) = e^{-h^2/2}$ be the characteristic function of the standard Gaussian distribution. Following is the main lemma, which will be used to prove Theorem~\ref{th1_speed_of_conv}.
\begin{lemma}
  \label{lemma1_growing_arg}
  Let~$f$ satisfy the assumptions of Theorem~\ref{th1_speed_of_conv}. Then, for every~$\varepsilon>0$ and~$\gamma \in [0,1/d]$ we have
  \begin{equation}
    \label{lemma1_growing_arg_form}
    \sup_{n}\sup_{|h|<\varepsilon n^\gamma} \left(n^{1-(d-1)\gamma} \left| \frac{\varphi_{f,n}(h)-\varphi_{\mathcal{N}}(h)}{h\varphi_{\mathcal{N}}(h)} \right|\right) < +\infty,
  \end{equation}
  where~$d=5$ for GUE and LUE, $d=3$ for JUE.
\end{lemma}

It is worth mentioning that an equivalent expression for~\eqref{variance_kf} is given by
\begin{equation}
  \label{quadr_series}
  K[f] = \frac{1}{4}\sum \limits_{j=1}^{\infty} j a_j^2,\quad a_j = \frac{2}{\pi} \int \limits_0^{\pi} f(\cos{s}) \cos{js}\, ds,
\end{equation}
from which~$K[f] \ge 0$ follows immediately. The~$a_j$ are the generalized Fourier coefficients with respect to the orthogonal system of the Chebyshev polynomials of the first kind~$\{T_j(x)\}_{j=0}^\infty$,
\begin{equation}
  f(x) = \sum \limits_{j=0}^{\infty} a_j T_j(x), \quad x \in [-1,1].
\end{equation}

Of special interest is the situation in which~$\gamma = 0$. The condition~$f(x) = O(e^{A\, V(x)})$ as~$x \to \infty$, $x \in \mathcal{I}$, turns out to be unnecessary, and we have another lemma.
\begin{lemma}
  \label{lemma2_fixed_h}
  Let~$f:\mathcal{I} \to \mathbb R$ be a locally H\"{o}lder continuous function admitting an analytic continuation into a complex neighborhood of the interval~$[-1,1]$. Then, for every~$\varepsilon>0$
  \begin{equation}
    \label{lemma2_fixed_h_form}
    \sup_{n}\sup_{|h|<\varepsilon} \left(n\left| \varphi_{f,n}(h)-\varphi_{\mathcal{N}}(h) \right|\right) < +\infty.
  \end{equation}
\end{lemma}
This lemma and the continuity theorem for characteristic functions (e.g., see~\cite{Feller2}) yield CLT for all three ensembles straightaway. On the other hand, CLT can also be derived from the asymptotics of the \textit{real} exponential moments~$\mexp[n]{e^{h S_{f,n}}}$, $h \in \mathbb{R}$, which is a special case of the results of Charlier and Gharakhloo~\cite{Charlier2019}. That said, we stress that the conclusion of Lemma~\ref{lemma1_growing_arg} does not follow from those results. Moreover, we do not see how to estimate the rate of convergence from asymptotics of the real exponential moments alone. This is why Lemma~\ref{lemma1_growing_arg} is essential to proving Theorem~\ref{th1_speed_of_conv}.

We also emphasize that since the existence of the exponential moments is not required, Lemma~\ref{lemma2_fixed_h} holds for a broader class of the test functions in the case of GUE and LUE, unlike the results from~\cite{Charlier2019}. For instance, let $T_k$ be the Chebyshev polynomial of the first kind and of degree~$k$. 
For~$k \ge 1$ set
\begin{equation}
  \varkappa_k = \left\{
    \begin{aligned}
      -&\delta_{k,2}/2  &&\mbox{for GUE},\\
      -&\delta_{k,1}/2 &&\mbox{for LUE},\\
      0& &&\mbox{for JUE},
    \end{aligned}
  \right.\quad
  \mu_k = \left\{
    \begin{aligned}
      &0  &&\mbox{for GUE},\\
      &(-1)^{k-1}\alpha/2 &&\mbox{for LUE},\\
      &((-1)^{k-1} \alpha -\beta)/2 &&\mbox{for JUE}.
    \end{aligned}
  \right.
\end{equation}
Introduce the diagonal  matrix~$\Sigma =  \frac{1}{4}\mathrm{diag}\{1,\ldots,l\}$ and the corresponding centered Gaussian distribution $N(0,\Sigma)$. The corollary below is straightforward to prove by using Lemma~\ref{lemma2_fixed_h}.
\begin{corollary}
  \label{th2}
  Set~$Y_k = \Tr{T_k(M)} -n \varkappa_k - \mu_k$. Then, the convergence in distribution holds
  \begin{equation}
    (Y_1, \ldots, Y_l) \overset{d}{\longrightarrow} N(0,\Sigma),\quad n\to \infty.
  \end{equation}
\end{corollary}

\section{Outline of proof and discussion}
\label{subsect:outl}
We mentioned earlier that our proof of Theorem~\ref{th1_speed_of_conv} uses Feller's smoothing inequality, which, we remind, gives a bound on the Kolmogorov--Smirnov distance between two distributions in terms of a certain integral involving the corresponding characteristic functions (for details, see~\cite{Feller2}). Naturally, the more one knows about the characteristic functions of the given distributions, the better bounds on the Kolmogorov--Smirnov distance one can get. To establish the results in Theorem~\ref{th1_speed_of_conv} we need to control the characteristic function~$\varphi_{f,n}(h)$ for small~$h$ (i.e., $h = O(1)$ as~$n \to \infty$) and for large~$h$ (i.e., $h = O(n^\gamma)$, $\gamma \in (0,1)$, as~$n \to \infty$), which is done by means of Lemma~\ref{lemma2_fixed_h} and~\ref{lemma1_growing_arg}, respectively. The final result is obtained by optimizing with respect to~$\gamma$ (see Section~\ref{sect_proof_of_th1} for details).

The proof of Lemmas~\ref{lemma1_growing_arg} and~\ref{lemma2_fixed_h} is based on asymptotic analysis of~$\varphi_{f,n}$ as~$n \to \infty$ via the well-known method called the Riemann--Hilbert (RH) approach.  The usual version of this approach includes the following steps. The first step is to consider a special \textit{deformation} of the original test function~$f$, in other words, a one-parametric family of functions~$\{f_t\}_{t \in [0,1]}$ such that~$f_0=0$ and~$f_1=f$. This deformation induces the deformation~$\{\varphi_{f_t,n}\}_{t \in [0,1]}$ of the characteristic function~$\varphi_{f,n}$. The next step is to establish a connection between~$\{\varphi_{f_t,n}\}_{t \in [0,1]}$ and the Hankel determinants, which is done via Andr\'{e}ief's identity (e.g., see~\cite{Andreev1886}). Then the further step is to connect the Hankel determinants with orthogonal polynomials and to form a special~$2 \times 2$ matrix out of these polynomials. This matrix turns out to satisfy a certain RH problem in the complex-analytic sense (e.g., see~\cite{Fokas1992}). The final step is to perform the asymptotic analysis of this problem (also known as the nonlinear steepest descent analysis of Deift and Zhou~\cite{Deift1993}) and to recover the asymptotics of the characteristic function~$\varphi_{f,n}$ by integrating (with respect to~$t \in [0,1]$) the special differential identity for~$\{\varphi_{f_t,n}\}_{t \in [0,1]}$. This identity yields the asymptotics of~$\log{\varphi_{f_1,n}} - \log{\varphi_{f_0,n}} = \log{\varphi_{f,n}}$.

All steps of the RH approach are fairly well-established in the literature, e.g., see~\cite{Beresticki2017, Charlier2018, Charlier2019, Claeys2011, Claeys2015,  Deift2011, Deift2014, Deift1999, Deift1999a, Fyodorov2016, Krasovsky2007, Kuijlaars2004, Vanlessen2005, Xu2014, Zhao2013}. It is worth mentioning, however, that the last step of the method is by far the most involved one. This step itself consists of the several substeps, equivalent transformations of the RH problem, during which the so-called global and local parametrices are constructed and the deformation of the contour is carried out. For the sake of the reader's convenience, we overview the RH approach and introduce a convenient notation in Section~\ref{RH_analysis}: this information is necessary to understanding Sections~\ref{lemma1_proof} and~\ref{sec_proof_lemma_gr_arg}.

One of the specific problems we have to address in our study is the presence of the imaginary exponent in the definition of~$\varphi_{f,n}$. This translates into the issue that the symbols of the Hankel determinants corresponding to~$\{\varphi_{f_t,n}\}_{t \in [0,1]}$ can have zeros in any fixed complex neighborhood of~$[-1,1]$, and thus the usual nonlinear steepest descent analysis does not go through (see discussion in Section~\ref{lemma1_proof}). To overcome this issue we employ an idea used earlier by Deift, Its, and Krasovsky in~\cite{Deift2014} while  studying Toeplitz determinants. Instead of only one deformation, $\{f_t\}_{t \in [0,1]}$, we construct a family (a ``chain") of deformations~$\left\{\{f_{l,t}\}_{t \in [0,1]}\right\}_{l=1,\ldots,q}$ --- by definition~$f_{l-1,1}=f_{l,0}$, $f_{0,0}:=1$, and $f_{q,1}:=f$ --- such that each individual~$f_{l,t}$ does not have zeros in a small enough complex neighborhood of~$[-1,1]$ for all~$t \in [0,1]$. Consequently, one can apply the usual RH analysis to each~$\{f_{l,t}\}_{t \in [0,1]}$ and recover the desired asymptotics of~$\varphi_{f,n}$ step by step, integrating the differential identity~$q$ times. Each integration gives the asymptotics of~$\log{\varphi_{f_{l,1},n}}-\log{\varphi_{f_{l,0},n}}$. By summing over all~$l=1,\ldots,q$, the asymptotic of~$\log{\varphi_{f_{q,1},n}}-\log{\varphi_{f_{0,0},n}} = \log{\varphi_{f,n}}$ is recovered. We emphasize that~$q$ is independent of~$n$ here, so the summation does not affect the order of the resulting error in the asymptotics. For details, see Section~\ref{lemma1_proof}.

Another issue that we face is connected with the fact that in order to obtain a bound on the rate of convergence, we need the asymptotics of the characteristic function~$\varphi_{f,n}(h)$ to be uniform in~$h$ for~$|h|<\varepsilon n^\gamma$, as~$n \to \infty$. In other words, we need a uniform asymptotic expansion of~$\varphi_{f,n}(n^\gamma h)$ for~$|h|<\varepsilon$. To our best knowledge, similar problems have not been addressed in the literature to date, and our idea is to construct a special deformation of~$\varphi_{f,n}(n^\gamma h)$ similar to that described in the previous paragraph. The difference is, however, that we allow~$l$ to range over all natural numbers, retaining the good control over the error in the corresponding asymptotics (for details, see Section~\ref{sec_proof_lemma_gr_arg}). Then it turns out that we can recover the desired asymptotics by carrying out the usual steepest descent analysis.

We noted (see Remark~\ref{rem_dep_on_loc_param}) that the bound on the rate of convergence in Theorem~\ref{th1_speed_of_conv} is due to the local behavior of the eigenvalues of M near the (asymptotic) edges of the spectrum, $x = \pm 1$. Indeed, it is well-known that the edge behavior of the eigenvalues can be given in terms of the local parametrices. And we will see in Section~\ref{sec_proof_lemma_gr_arg} that these local parametrices is exactly what determine the error term in the asymptotics of the characteristic function, and thus the bound on the rate of convergence (see Remark~\ref{remark: loc_par_define_rate}). Depending on the type of an edge, soft or hard, the local parametrices differ. And it turns out that if at least one soft edge is present, the bound on the rate is of lower order than otherwise.

It is worth mentioning that a somewhat similar effect has been observed in~\cite[Remark~1.3]{Lambert2019} while studying the quadratic Kantorovich distance between the standard Gaussian law and the law of a linear statistic for~$\beta$-ensembles with one-cut potentials. It has been shown there that when a test function is supported away from the edges, the rate of convergence becomes higher than otherwise. In the particular case of GUE and polynomial test functions, the optimal rate of~$O(n^{-1})$ was obtained. Based on this information, we expect that the bounds we have are not sharp; however, we have not been able to establish this rigorously.

We stress that the bound we obtain is of order~$O(n^{-1/d})$ and not~$O(n^{-1})$, as one would expect, because of the two factors. First, an essential step of our approach is to shrink the contours in the local RH analysis (see~\eqref{eq_our_proof_unb_shrink} and below), which causes the approximation error for the local parametrices to be of greater order than~$O(n^{-1})$. Second, we have to add up a large number of terms --- the number of terms is of order~$O(n^\gamma)$ --- to obtain the final asymptotics (see~\eqref{eq_148} and below). This gives an extra factor of~$n^\gamma$ in front of the resulting error term.  Both steps seem to be imperative, and we do not see how to avoid them.

At the end of this section we mention a few works related to the asymptotic analysis of the Hankel determinants. Johansson~\cite{Johansson1998} deals with the matrix models with the continuous weight on the whole real line. Vanlessen~\cite{Vanlessen2005} studies the Plancherel--Rotach asymptotics for the orthogonal polynomials with a Laguerre-type weight. Zhao, Cao and Dai~\cite{Zhao2013} obtained the asymptotic expansion of the partition function of a Laguerre-type model. The case of a Laguerre-type singularly perturbed weight was studied by Xu, Dai and Zhao in~\cite{Xu2014}, where the connection was found between the Painlev\'{e} III transcendent and the behavior of the leading and recurrence coefficients of the corresponding orthogonal polynomials. Lyu and Chen~\cite{Lyu2017} studied the distribution of the largest eigenvalue in the Laguerre unitary ensembles. In the physical literature, the connection between eigenvalue statistics and wireless relaying has been studied by Chen and Lawrence~\cite{Chen1998},  Chen, Haq and McKay~\cite{Chen2013}.

\section{Analysis of the Riemann--Hilbert problem}
\label{RH_analysis}

Here, we give an overview of the classical RH approach (also, see references in Section~\ref{subsect:outl}), which is an essential building block of our proofs. Let~$\tilde{f}$ be a complex-valued function on $\mathcal{I}$ such that the following holds.
\begin{assumption}
  \label{assump:A1}
  The function~$\exp{\tilde{f}}$ is  locally H\"{o}lder continuous on~$\mathcal{I}$. Additionally, for GUE and LUE $\tilde{f}$ satisfies~$\max\!{\{\Re{\tilde{f}(x)},0\}} = O(V(x))$ as~$|x| \to +\infty$, $x \in \mathcal{I}$.
\end{assumption}
\begin{assumption}
  \label{assump:A3}
  The function $\tilde{f}$~admits an analytic continuation into a complex neighborhood of~$[-1,1]$.
\end{assumption}

First, we connect the expectation~$\mexp[n]{e^{\Tr{\tilde{f}(M)}}}$ and the Hankel determinants. Passing to the radial part in~\eqref{nu-measure}, we have
\begin{equation}
  \label{math-exp}
  \mexp[n]{e^{\Tr{\tilde{f}(M)}}} = \frac{1}{Z_{n}} \int \limits_{\mathcal{I}^n} e^{\sum \limits_j (\tilde{f}(\lambda_j) - Q_n(\lambda_j))} \prod \limits_{j<k} (\lambda_k -\lambda_j)^2\, d\lambda_1\cdot \ldots \cdot d\lambda_n;
\end{equation}
and then Andr\'{e}ief's identity (e.g., see~\cite{Andreev1886}) gives
\begin{equation}
  \label{exp-via-hankel}
  \mexp[n]{e^{\Tr{\tilde{f}(M)}}} = \frac{H_{n,n}[\tilde{f}]}{H_{n,n}[0]},
\end{equation}
where~$H_{n,m}[\tilde{f}] = \det\{\mu_{j+k-2}^{(m)}\}_{j,k=1}^n$ is the Hankel determinant with the symbol given by
\begin{equation}
  \label{weight}
  w_m(x) = e^{\tilde{f}(x)- Q_m(x)}, \quad x \in \mathcal{I},
\end{equation}
and the~$\mu_j^{(m)}$ are the moments
\begin{equation}
  \label{moments}
  \mu_j^{(m)} = \int\limits_{\mathcal{I}} x^j w_m(x)\, dx.
\end{equation}
Assumption~\ref{assump:A1} guarantees that the integral in~\eqref{moments} exists and the transition from~\eqref{math-exp} to~\eqref{exp-via-hankel} is legitimate. Assumption~\ref{assump:A3} is needed while deforming the contour (and the RH problem itself) in Section~\ref{subset_sec_transf} and while constructing the local parametrices in Sections~\ref{subsect_locp1} and~\ref{subsect_locm1}.

Now, assume additionally that the following holds.
\begin{assumption}
  \label{assump:A4}
  $H_{n-1,n}[\tilde{f}]\ne 0$ and~$H_{n,n}[\tilde{f}] \ne 0$  for all~$n$ large enough. 
\end{assumption}
This implies that the (monic) orthogonal polynomials~$\pi^{(k)}_{n}(x) = x^k+\ldots$ with respect to the weight~\eqref{weight} are well-defined for~$k=n-1, n$, where~$n$ is large enough.

Consider the matrix function
\begin{equation}
  \label{eq_probly_solution}
  Y_{n}(z) = \begin{bmatrix}
    \pi^{(n)}_{n}(z) & \mathcal{C}(\pi^{(n)}_{n} w_{n})(z)\\
    \beta_{n-1,n}\pi^{(n-1)}_{n}(z) &  \beta_{n-1,n} \mathcal{C}(\pi^{(n-1)}_{n} w_{n})(z)
  \end{bmatrix}, \quad z \in \mathbb{C} \setminus \mathcal{I},
\end{equation}
where~$\beta_{n,m} = - 2\pi i \gamma_{n,m}^2$, $\gamma_{n,m}^2 =\frac{H_{n,m}[\tilde{f}]}{H_{n+1,m}[\tilde{f}]}$, and~$\mathcal{C}$ is the Cauchy-type integral
\begin{equation}
  \mathcal{C}(g)(z) = \frac{1}{2 \pi i} \int\limits_{\mathcal{I}} \frac{g(s)}{s-z} \, ds, \quad z \in \mathbb{C} \setminus \mathcal{I}.
\end{equation}
Further, we will often drop the subscript~$n$ to make the notation lighter.

Also let~$\mathring{\mathcal{I}}$ be the set of the interior points of~$\mathcal{I}$. Because of Assumption~\ref{assump:A1}, the upper and lower limits~$Y^\pm(x) = \lim\limits_{z \to x \pm i0} Y(z)$, $x \in \mathring{\mathcal{I}}$, are well defined pointwise (e.g., see~\cite{Gakhov_book}). Then it is straightforward to check that~$Y(z)$ solves the following RH problem.
\begin{problem}[Y-RH]
  \namedlabel{problem_y-rh}{Problem Y-RH}
  \leavevmode
  \begin{enumerate}[label=\textnormal{({\roman*})},ref=Y-RH-{\roman*}]
  \item  \label{yrh_cond1} $Y(z)$ is analytic in~$\mathbb{C} \setminus \mathcal{I}$;
  \item  \label{yrh_cond2} $Y^+(x) = Y^-(x) J_Y(x),\, x \in \mathring{\mathcal{I}}$, where~$J_Y(x) =
    \begin{bmatrix}
      1 & w_n(x)\\
      0&1
    \end{bmatrix}$;
  \item  \label{yrh_cond3} $Y(z) = (I + \mathcal{O}(1/z))z^{n \sigma_3}$ as~$z \to \infty$, where~$\sigma_3 = \begin{bmatrix}
      1 & 0\\
      0& -1
    \end{bmatrix}$;
  \item  \label{yrh_cond4}
    \begin{equation}
      Y(z) = \left\{
        \begin{aligned}
          &\mathcal{O}(1)&&\mbox{ for GUE},\\
          &\begin{bmatrix}
            O(1) & O(1) + O(|z+1|^\alpha)\\
            O(1) & O(1) + O(|z+1|^\alpha)
          \end{bmatrix}, &\alpha \ne 0, &\mbox{ for LUE, JUE,}\\
          &\begin{bmatrix}
            O(1) & O(\log{|z+1|})\\
            O(1) & O(\log{|z+1|})
          \end{bmatrix}, &\alpha =0, &\mbox{ for LUE, JUE},
      \end{aligned}
    \right.
  \end{equation}
    as $z \to -1, z \in \mathbb{C} \setminus \mathcal{I}$,
    \begin{equation}
      Y(z) = \left\{
        \begin{aligned}
          &\mathcal{O}(1) &&\mbox{ for GUE, LUE},\\
          &\begin{bmatrix}
            O(1) & O(1) + O(|z-1|^\beta)\\
            O(1) & O(1) + O(|z-1|^\beta)
          \end{bmatrix}, &\beta \ne 0,&\mbox{ for JUE,}\\
          &\begin{bmatrix}
            O(1) & O(\log{|z-1|})\\
            O(1) & O(\log{|z-1|})
          \end{bmatrix}, &\beta =0, &\mbox{ for JUE,}\\
        \end{aligned}
      \right.
    \end{equation}
    as $z \to 1, z \in \mathbb{C} \setminus \mathcal{I}$.
  \end{enumerate}
\end{problem}
From here on out, for the sake of brevity we use the notation
\begin{equation}
  \mathcal{O}(1/z) = \begin{bmatrix}
    O(1/z) & O(1/z)\\
    O(1/z)& O(1/z)
  \end{bmatrix}.
\end{equation}
Notice that the~$O$ terms above can depend on~$n$. This does not cause any trouble, however, because in our further analysis we use~\eqref{yrh_cond3} and~\eqref{yrh_cond4} only when~$n$ is fixed.

We also mention that, essentially, a link between orthogonal polynomials and RH problems was established by Fokas, Its, and Kitaev~\cite{Fokas1992} while studying the Hermitian matrix model for 2D quantum gravity.
\begin{remark}
  \label{remark:no_assump_A4}
  For our further analysis Assumption~\ref{assump:A4} is not needed. Above, we used it exclusively to write~\eqref{eq_probly_solution}.
\end{remark}

It is well-known that~\ref{problem_y-rh} has a unique solution, which satisfies~$\det{Y(z)}=1$ (e.g., see~\cite[p.~44]{Deift_book}). Indeed, it follows from~\eqref{yrh_cond1} that~$\det{Y(z)}$ is an analytic function in~$\mathbb{C} \setminus \mathcal{I}$; moreover due to~\eqref{yrh_cond2}, $\det{Y(z)}$ has no jumps over~$\mathring{\mathcal{I}}$. For LUE and JUE, it is possible, though, that~$\det{Y(z)}$ has isolated singularities at~$z=\pm 1$; nevertheless, these singularities are removable because of~\eqref{yrh_cond4}. So, $\det{Y(z)}$ turns out to be an entire function in all three cases. Finally, \eqref{yrh_cond3} and Liouville's theorem ensure that~$\det{Y(z)}=1$ for all~$z \in \mathbb{C}$. In particular, this means that~$Y(z)$ is invertible, i.e., $(Y(z))^{-1}$ is well-defined. Now, suppose that there are two solutions of~\ref{problem_y-rh}, $Y_1(z)$ and~$Y_2(z)$. Using  the similar reasoning as above, we see that~$Y_1(z)(Y_2(z))^{-1}$ is the identity matrix, and thus the solution is unique.

In the sections that follow, we describe the steepest descent analysis of~\ref{problem_y-rh}, which, we remind, involves a series of equivalent transformations of this problem, from one RH problem to another. Our ultimate goal is to have a problem normalized at infinity with the corresponding jump matrix converging to the identity matrix uniformly on the contour as~$n \to \infty$. In this case the well-known results of the theory of small-norm RH problems can be applied (e.g., see~\cite{Deift_book, Its2011}). We start by normalizing~\ref{problem_y-rh} at infinity.

\subsection{First transformation:  normalization at~$z = \infty$}
\label{subsect_eq_measure}
Consider the equilibrium measure~$\nu(dx)$ corresponding to a potential~$V(x)$. This measure is the unique solution to the variational problem
\begin{equation}
  \label{equilib_mes}
  \iint \limits_{x\ne y} \log{\frac{1}{|x-y|}} \mu(dx)\mu(dy) + \int V(x) \mu(dx) \to \min,
\end{equation}
where the minimization is done over the (convex) set of the probability measures~$\mu(dx)$ supported on~$\Sigma \subset \mathcal{I}$. The optimality conditions, following from the corresponding variational inequality, can be written as
\begin{equation}
  \label{lagr_cond}
  \begin{aligned}
    &2 \int \log{\frac{1}{|x-y|}}\, \mu(dy) + V(x) = l_R, \quad x \in \Sigma,\\
    &2 \int \log{\frac{1}{|x-y|}}\, \mu(dy) + V(x) \ge l_R, \quad x \in \mathcal{I} \setminus \Sigma,
  \end{aligned}
\end{equation}
where~$l_R$ is a real number, called the modified Robin constant (see~\cite{Saff&Totik}).

The latter system can be solved explicitly for each of the three potentials in~\eqref{V_omega} (e.g., see~\cite{Charlier2019, Deift1999a}). The corresponding equilibrium measures are supported on~$[-1,1]$. And it turns out that they are absolutely continuous with respect to the Lebesgue measure, that is, $\mu(dx) = \psi(x)\ dx$. The density~$\psi(x)$, $x \in [-1,1]$, and the corresponding Robin constant~$l_R$ are
\begin{equation}
  \label{eq:equil_dens}
  \psi(x) = \left\{
    \begin{aligned}
      &\frac{2}{\pi} \sqrt{1-x^2} &&\mbox{for GUE},\\
      &\frac{1}{\pi} \sqrt{\frac{1-x}{1+x}} &&\mbox{for LUE},\\
      &\frac{1}{\pi} \frac{1}{\sqrt{1-x^2}} &&\mbox{for JUE},
    \end{aligned}
  \right.\quad
  l_R = \left\{
    \begin{aligned}
      &1 + 2 \log{2} &&\mbox{for GUE},\\
      &2 + 2 \log{2} &&\mbox{for LUE},\\
      &2 \log{2} &&\mbox{for JUE}.
    \end{aligned}
  \right.
\end{equation}

It is worth mentioning that in all three cases the~$\mu$'s are, in fact, probability distributions; in the respective order, Wigner's semicircle distribution, the Marchenko--Pastur distribution, and the arcsine distribution. Each of them is the weak limit of the normalized counting measure of eigenvalues for the corresponding random matrices (e.g., see~\cite{Deift_book, Pastur_book}). This explains why~$\varkappa[\cdot]$ defined in~\eqref{eq:ch1_leading_asymp_kappa}, the expectation with respect to the equilibrium measure, is naturally present in~\eqref{eq_eq_Sfn} and why the leading term of the asymptotics of~$\mexp[n]{\Tr{f(M)}}$ as~$n \to \infty$ is~$n \varkappa[f]$.

Now, consider the logarithmic potential
\begin{equation}
  \label{log_potential}
  g(z) = \int\limits_{-1}^1 \log(z-s) \psi(s)\ ds, \quad z \in \mathbb{C}\setminus (-\infty,1],
\end{equation}
and the auxiliary function
\begin{equation}
  \label{phi}
  \phi(z) = \left\{
    \begin{aligned}
      &4 \int \limits_1^z \sqrt{s^2-1}\ ds &&\mbox{for GUE},\\
      &2 \int \limits_1^z \sqrt{\frac{s-1}{s+1}}\ ds &&\mbox{for LUE},\\
      &-2 \int \limits_1^z \frac{1}{\sqrt{s^2-1}}\ ds &&\mbox{for JUE},
    \end{aligned}
  \right.\quad
  z \in \mathbb{C}\setminus (-\infty,1],
\end{equation}
where the principal branches of the logarithm and roots are used. By appealing to these definitions and using~\eqref{lagr_cond}, one can easily show that~$g(z)$ and~$\phi(z)$ are analytic in~$\mathbb{C}\setminus (-\infty,1]$; moreover, the following identities hold:
\begin{equation}
  \label{g_func_prop}
  \begin{aligned}
    &2g(z)- V(x) + l_R = -\phi(z), &&z \in \mathbb{C}\setminus (-\infty,1],\\
    &g^+(x)+g^-(x)- V(x) + l_R = 0, &&x \in [-1,1],\\
    &g^+(x)-g^-(x) = -\phi^+(x) = \phi^-(x), &&x \in [-1,1],\\
    &g^+(x)-g^-(x) = 0, &&x \in [1,+\infty),\\
    &g^+(x)-g^-(x) = 2 \pi i, &&x \in (-\infty,1],
  \end{aligned}
\end{equation}
where the superscript~$+$ and~$-$ are used to denote the upper and lower half-plane limits, respectively.

Now, we are ready to carry out the first step of the steepest descent analysis, the change of variables in~\ref{problem_y-rh}:
\begin{equation}
  \label{U-trans}
  U(z) = e^{n l_R \sigma_3/2}Y(z)e^{-n(l_R/2 + g(z)) \sigma_3}.
\end{equation}
Clearly, $U(z)$ is analytic in~$\mathbb{C} \setminus \mathcal{I}$ thanks to~\eqref{g_func_prop}. Besides, since~$g(z) = \log(z) + O(1/z)$ as~$z \to \infty$, which follows from~\eqref{log_potential}, one also has~$U(z) = I +O(1/z)$ as~$z \to \infty$. So, $U(z)$ is normalized at the infinity.

For convenience, introduce
\begin{equation}
  \label{chi}
  \chi(x) = e^{\tilde{f}(x) + \omega(x)}.
\end{equation}
Then it is easy to check directly that~$U(z)$ is the solution of the following RH problem.
\begin{problem}[U-RH]
  \namedlabel{problem_u-rh}{Problem~U-RH}
  \leavevmode
  \begin{enumerate}[label=\textnormal{({\roman*})},ref=U-RH-{\roman*}]
  \item $U(z)$ is analytic in~$\mathbb{C} \setminus \mathcal{I}$;
  \item $U^+(x) = U^-(x) J_U(x),\, x \in \mathring{\mathcal{I}}$,
    \begin{equation}
      \label{jump_U}
      J_U(x) =\left\{
        \begin{aligned}
          &\begin{bmatrix}
            1 & \chi(x) e^{-n\phi(x)}\\
            0 & 1
          \end{bmatrix}, &&x \in \mathcal{I}\setminus [-1, 1],\\
          &\begin{bmatrix}
            e^{n\phi^+(x)} & \chi(x)\\
            0 & e^{n\phi^-(x)}
          \end{bmatrix}, &&x\in (-1,1);
        \end{aligned}
      \right.
    \end{equation}
  \item $U(z) = (I + \mathcal{O}(1/z))$ as~$z \to \infty$;
  \item The behavior of~$U(z)$ as~$z \to \pm 1$ is the same as that of~$Y(z)$ in \textnormal{~\ref{problem_y-rh}}.
  \end{enumerate}
\end{problem}
The formula~\eqref{jump_U} for the jump matrix~$J_U(x)$ is easily obtained by using~\eqref{yrh_cond2} and~\eqref{g_func_prop}:
\begin{equation}
  \begin{aligned}
    J_U(x) &= (U^-(x))^{-1}U^+(x) = e^{n(l_R/2 + g^-(x)) \sigma_3} J_Y(x) e^{-n(l_R/2 + g^+(x)) \sigma_3}\\
    &= 
    \begin{bmatrix}
      e^{-n(g^+(x)-g^-(x))} & \chi(x) e^{n(g^+(x)+g^-(x)- V(x)+l_R)}\\
      0&e^{n(g^+(x)-g^-(x))}
    \end{bmatrix}\\
    &=\left\{
      \begin{aligned}
        &\begin{bmatrix}
          1 & \chi(x) e^{-n\phi(x)}\\
          0 & 1
        \end{bmatrix}, &&x \in \mathcal{I}\setminus [-1, 1],\\
        &\begin{bmatrix}
          e^{n\phi^+(x)} & \chi(x)\\
          0 & e^{n\phi^-(x)}
        \end{bmatrix}, &&x\in (-1,1).
      \end{aligned}
    \right.
  \end{aligned}
\end{equation}

We highlight that for GUE and LUE, $\phi(z)$ is defined in such a way that~$\phi(x)>0$ for~$x \in \mathcal{I}\setminus [-1, 1]$; therefore, from the definition~\eqref{jump_U} it is immediate that~$J_U(x) \to I$, $x \in \mathcal{I}\setminus [-1, 1]$, as~$n \to \infty$. The rate of this convergence is exponential, but the convergence itself is not uniform on~$\mathcal{I}\setminus [-1, 1]$. Also, for all three ensembles, $\phi^\pm(x)$ are imaginary if~$x \in (-1,1)$, and thus~$J_U(x)$ exhibit oscillatory behavior as~$n \to \infty$ and does not converge anywhere. This means that the theory of the small-norm RH problems cannot be applied, which is why the next step of the RH approach is needed.

\subsection{Second transformation:  deformation of the RH problem}
\label{subset_sec_transf}
Consider the jump matrix~\eqref{jump_U} on~$(-1,1)$, and observe that a simple matrix identity takes place:
\begin{equation}
  \label{matr_fact}
  \begin{aligned}
    \begin{bmatrix}
      e^{n \phi^+(x)} & \chi(x)\\
      0&e^{n\phi^-(x)}
    \end{bmatrix} &=
    \begin{bmatrix}
      1 & 0\\
      \frac{e^{n \phi^-(x)}}{\chi(x)} &1
    \end{bmatrix}
    \begin{bmatrix}
      0 & \chi(x)\\
      -\frac{1}{\chi(x)}&0
    \end{bmatrix}
    \begin{bmatrix}
      1 & 0\\
      -\frac{e^{n \phi^+(x)}}{\chi(x)} &1
    \end{bmatrix}^{-1} \\
    &=: J_T^-(x) J_T^{o}(x)J_T^+(x).
  \end{aligned}
\end{equation}

\begin{figure}[h!]
  \centering
  \begin{subfigure}[b]{0.32\textwidth}
    \centering
    \begin{tikzpicture}[thick,scale=0.8]
      \begin{scope}[decoration={markings,mark= at position 0.5 with {\arrow{stealth}}}]
        \draw[postaction={decorate}] (-1,0)-- (0,0);
        \draw[postaction={decorate}] (0,0) -- (3,0);
        \draw[postaction={decorate}] (3,0)-- (4,0);
        \draw[postaction={decorate}] (0,0) to [out=50, in=130] (3,0);
        \draw[postaction={decorate}] (0,0) to [out=-50, in=-130] (3,0);
        \node at (3,-0.3) {$1$};
        \node at (-0.2,-0.3) {$-1$};
        \node at (1.6,-0.35) {$\Omega^-$};
        \node at (1.6,0.35) {$\Omega^+$};
        \node at (2.5,0.8) {$L_+$};
        \node at (2.5,-0.8) {$L_-$};
      \end{scope}
    \end{tikzpicture}
    \caption{GUE}
  \end{subfigure}
  \begin{subfigure}[b]{0.32\textwidth}
    \centering
    \begin{tikzpicture}[thick,scale=0.8]
      \begin{scope}[decoration={markings,mark= at position 0.5 with {\arrow{stealth}}}]
        \draw[draw=none] (-1,0)-- (0,0);
        \draw[postaction={decorate}] (0,0) -- (3,0);
        \draw[postaction={decorate}] (3,0)-- (4,0);
        \draw[postaction={decorate}] (0,0) to [out=50, in=130] (3,0);
        \draw[postaction={decorate}] (0,0) to [out=-50, in=-130] (3,0);
        \node at (3,-0.3) {$1$};
        \node at (-0.2,-0.3) {$-1$};
        \node at (1.6,-0.35) {$\Omega^-$};
        \node at (1.6,0.35) {$\Omega^+$};
        \node at (2.5,0.8) {$L_+$};
        \node at (2.5,-0.8) {$L_-$};
      \end{scope}
    \end{tikzpicture}
    \caption{LUE}
  \end{subfigure}
  \begin{subfigure}[b]{0.32\textwidth}
    \centering
    \begin{tikzpicture}[thick,scale=0.8]
      \begin{scope}[decoration={markings,mark= at position 0.5 with {\arrow{stealth}}}]
        \draw[postaction={decorate}] (0,0) -- (3,0);
        \draw[postaction={decorate}] (0,0) to [out=50, in=130] (3,0);
        \draw[postaction={decorate}] (0,0) to [out=-50, in=-130] (3,0);
        \node at (3,-0.3) {$1$};
        \node at (-0.2,-0.3) {$-1$};
        \node at (1.6,-0.35) {$\Omega^-$};
        \node at (1.6,0.35) {$\Omega^+$};
        \node at (2.5,0.8) {$L_+$};
        \node at (2.5,-0.8) {$L_-$};
      \end{scope}
    \end{tikzpicture}
    \caption{JUE}
  \end{subfigure}
  \caption{The deformed contour~$L$ of~\ref{problem_t-rh}.}
  \label{T_cont}
\end{figure}
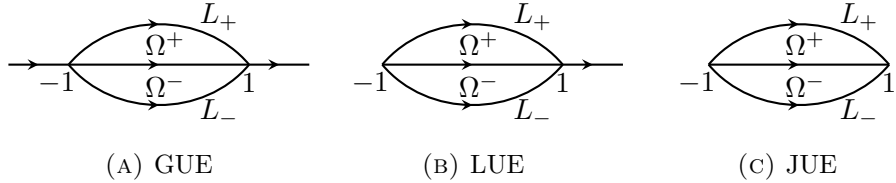

Then deform~$\mathcal{I}$ into the lens-shaped contour in Fig.~\ref{T_cont} (notice that the lips~$L_\pm$ do not include the edges~$\pm 1$) and set~$T(z)$ to be
\begin{equation}
  \label{T-def}
  T(z) = \left\{
    \begin{aligned}
      &U(z),&& z \in \Omega=\mathbb{C} \setminus (\mathcal{I} \cup \overline{\Omega^+ \cup \Omega^-}),\\
      &U(z) (J_T^+(z))^{-1}, && z \in \Omega^+,\\
      &U(z) J_T^-(z), && z \in \Omega^-.
    \end{aligned}
  \right.
\end{equation}
For~\eqref{T-def} to make sense, we need Assumption~\ref{assump:A3}, which guarantees the existence of an analytic continuation of~$\tilde{f}(x)$ into some complex neighborhood of~$[-1,1]$. Without loss of generality, the lens is embedded in this neighborhood, and we see from~\eqref{chi} that~$\chi(z)$ is a well-defined non-vanishing analytic function. Consequently, $J_T^\pm(x)$ also have analytic continuations~$J_T^\pm(z)$ from the interval~$(-1,1)$ into the respective domains~$\Omega^\pm$, moreover~$J_T^\pm(z)$ are continuous up to the boundary. We note that the analytic continuation of a function and the function itself are denoted by the same symbol; also, in order to define~$\chi(z)$ (see the formulas~\eqref{chi} and~\eqref{Omega_omega}), we use the principal branch of the logarithm on~$\mathbb{C} \setminus (-\infty,0]$.

It is readily verified that~$T(z)$ solves the following RH problem.
\begin{problem}[T-RH]
  \needspace{5ex}
  \namedlabel{problem_t-rh}{Problem~T-RH}
  \leavevmode
  \begin{enumerate}[label=\textnormal{({\roman*})},ref=T-RH-{\roman*}]
  \item \label{trh_cond1} $T(z)$ is analytic in~$\mathbb{C} \setminus L$;
  \item \label{trh_cond2} $T^+(z) = T^-(z) J_T(z),\, z \in L$,
    \begin{equation}
      \label{Y2RH}
      J_T(z) =\left\{
        \begin{aligned}
          &J_T^+(z),&& z \in L_+,\\
          &J_T^{o}(z),&& z \in (-1,1),\\
          &J_T^-(z),&& z \in L_-,\\
          &J_U(z), && z \in \mathcal{I} \setminus [-1,1];
        \end{aligned}
      \right.
    \end{equation}
  \item \label{trh_cond3} $T(z) = (I + \mathcal{O}(1/z))$ as~$z \to \infty$;
  \item \label{trh_cond4} If approaching from~$\Omega$, the behavior of~$T(z)$ in the neighborhoods of points~$\pm1$ is the same as that of~$U(z)$. If approaching from~$\Omega^+$ and~$\Omega^-$, the behavior can be obtained by multiplying by the corresponding jump matrix, as in~\eqref{T-def}.
  \end{enumerate}
\end{problem}
We notice that an informal interpretation of the transformation~\eqref{T-def}, in view of~\eqref{Y2RH}, is that it is meant to ``spread" the initial jump~\eqref{matr_fact}, originally over~$\mathcal{I}$ only, over the three contours~$L_+$, $L_-$, and~$\mathcal{I}$.

By the definition~\eqref{phi}, it is clear that~$\Re{\phi(z)}<0$ on~$L_\pm$. Consequently, $J_T(z) \to I$. Again, the rate of this convergence is exponential, but the convergence itself is not uniform on~$L_\pm$. That said, we note that extracting arbitrarily small neighborhoods~$\Omega_{\pm 1}$ of~$\pm1$ renders this convergence uniform on~$(L_+ \cup L_-) \setminus (\Omega_{1} \cup \Omega_{-1})$. We will construct the exact asymptotic solutions of~\ref{problem_t-rh} in~$\Omega_{\pm 1}$, the local parametrices, in Sections~\ref{subsect_locp1} and~\ref{subsect_locm1}.

Noticing that~$J_T(z) = J_U(z)$ if~$x \in \mathcal{I} \setminus [-1,1]$, we conclude that~$J_T(z) \to I$ uniformly on~$(L_+ \cup L_- \cup (\mathcal{I}\setminus [-1,1]))\setminus (\Omega_{1} \cup \Omega_{-1})$. Therefore, in view of the theory of small-norm RH problems, it makes sense to consider the limiting RH problem on the rest of the contour~$(-1,1)$ (we ignore~$\Omega_{1}$ and~$\Omega_{-1}$ since they can be arbitrarily small). The corresponding jump matrix is~$J_T^{o}(x)$. In the next section, the construction of the solution of this limiting problem, the so-called global parametrix, will be carried out.

\subsection{Global parametrix}
The construction of the global parametrix can be given explicitly. And this is one of the reasons for the success of the RH method. Consider the following problem.
\begin{problem}[N-RH]
  \namedlabel{problem_n-rh}{Problem~N-RH}
  \leavevmode
  \begin{enumerate}[label=\textnormal{({\roman*})},ref=N-RH-{\roman*}]
  \item \label{nrh_cond1} $N(z)$ is analytic in~$\mathbb{C} \setminus [-1,1]$;
  \item \label{nrh_cond2} $N^+(x) = N^-(x) J_T^{o}(x),\ x \in (-1,1)$,
  \item \label{nrh_cond3} $N(z) = (I + \mathcal{O}(1/z))$ as~$z \to \infty$.
  \end{enumerate}
\end{problem}
In order to find a solution of~\ref{problem_n-rh}, observe that
\begin{equation}
  J_T^{o}(x) = (\chi(x))^{\sigma_3} \cdot
  \begin{bmatrix}
    0 & 1\\
    -1 & 0
  \end{bmatrix},
\end{equation}
and introduce the Szeg\H{o} function~$D(z)$ by the formula
\begin{equation}
  \label{szego_func}
  D(z) = \exp{\left(\frac{1}{2 \pi} \sqrt{z^2-1} \int\limits_{-1}^1 \frac{\omega(s) + \tilde{f}(s)}{\sqrt{1-s^2}} \cdot \frac{ds}{z-s}\right)}, \quad z \in \mathbb{C}\setminus [-1,1],
\end{equation}
where~$\omega(x)$ is given by~\eqref{Omega_omega} and the principal branch of the root is used. Clearly, this function is analytic in~$\mathbb{C} \setminus [-1,1]$, and thanks to the Sokhotski--Plemelj formulas, one sees that
\begin{equation}
  \label{RH_mult_prob}
  D^+(x) D^-(x) = \chi(x), \quad x \in (-1,1).
\end{equation}

The expression~\eqref{szego_func} can be easily factorized~$D(z) = D_1(z) D_2(z)$, where
\begin{equation}
  \begin{aligned}
    \label{D12-funct}
    &D_1(z) = \exp{\left(\frac{1}{2 \pi} \sqrt{z^2-1} \int\limits_{-1}^1 \frac{\omega(s)}{\sqrt{1-s^2}} \cdot \frac{ds}{z-s} \right)},\\
    &D_2(z) = \exp{\left(\frac{1}{2 \pi} \sqrt{z^2-1} \int\limits_{-1}^1 \frac{\tilde{f}(s)}{\sqrt{1-s^2}} \cdot \frac{ds}{z-s} \right)}.
  \end{aligned}
\end{equation}
Moreover, it is known that~$D_1(z)$ has a simple explicit representation (e.g., see~\cite{Charlier2019, Vanlessen2005}). For the sake of completeness, we derive this representation here, using the residue theory. For definiteness, let us consider LUE, so
\begin{equation}
  \omega(x) = \alpha \log{(1+x)}.
\end{equation}
First, we introduce the change of variables~$s \mapsto \frac{1-s^2}{1+s^2}$. Then the integral in~$D_1(z)$ becomes
\begin{equation}
  I = \frac{\alpha}{2\pi}\int\limits_{-1}^1 \frac{\log{(1+s)}}{\sqrt{1-s^2}} \cdot \frac{ds}{z-s}  = -\frac{\alpha}{2\pi}\int \limits_{-\infty}^{+\infty} \frac{\log{((s^2+1)/2)}}{(z+1)s^2+z-1}\, ds.
\end{equation}
Calculating the corresponding residues, we obtain
\begin{equation}
  \begin{aligned}
    I &=  \frac{\alpha \log{(z+1)}}{2\sqrt{z^2-1}} - \alpha \int \limits_{1}^{+\infty} \frac{ds}{(z+1)s^2 -z+1} \\
    &= \frac{\alpha}{2\sqrt{z^2-1}} \left(\log{(z+1)} + \log{\frac{\sqrt{z+1}-\sqrt{z-1}}{\sqrt{z+1}+ \sqrt{z-1}}}\right),
  \end{aligned}
\end{equation}
where the principal branches of the logarithms and roots are used. Finally, we see that
\begin{equation}
  D_1(z) = (z+1)^{\alpha/2} \left(\frac{\sqrt{z+1}-\sqrt{z-1}}{\sqrt{z+1}+ \sqrt{z-1}}\right)^{\alpha/2}.
\end{equation}

It is interesting to point out that the function
\begin{equation}
  z \mapsto \frac{\sqrt{z+1}-\sqrt{z-1}}{\sqrt{z+1}+ \sqrt{z-1}}
\end{equation}
is a conformal map of the complex plane with the slit~$[-1,1]$ onto the interior of the standard unit disk without zero. Thereby, this function maps any loop going around~$[-1,1]$ counterclockwise into a loop inside the disk going around zero clockwise. In this way the jump of~$(z+1)^{\alpha/2}$ is compensated with the jump of~$\left(\frac{\sqrt{z+1}-\sqrt{z-1}}{\sqrt{z+1}+ \sqrt{z-1}}\right)^{\alpha/2}$, which makes~$D_1(z)$ analytic in~$\mathbb{C} \setminus [-1,1]$.

In the similar way one can handle the case of JUE. The final result is
\begin{equation}
  \label{glob_par_d1}
  D_1(z) = \left\{
    \begin{aligned}
      &1 &&\mbox{for GUE},\\
      &(z+1)^{\alpha/2} \left(\frac{\sqrt{z+1}-\sqrt{z-1}}{\sqrt{z+1}+ \sqrt{z-1}}\right)^{\alpha/2} &&\mbox{for LUE},\\
      &\begin{aligned}
        (z&+1)^{\alpha/2} (z-1)^{\beta/2}\\
        &\times \left(\frac{\sqrt{z+1}-\sqrt{z-1}}{\sqrt{z+1}+ \sqrt{z-1}}\right)^{(\alpha+\beta)/2}
      \end{aligned}&&\mbox{for JUE}.
    \end{aligned} 
  \right.
\end{equation}

Further, we will also need the formula for~$D(\infty)$. Clearly,
\begin{equation}
    \label{glob_param_D_infty}
  D(\infty) = D_1(\infty) D_2(\infty),
\end{equation}
where
\begin{equation}
  D_1(\infty) =\left\{
    \begin{aligned}
      &1 &&\mbox{for GUE},\\
      &2^{-\alpha/2} &&\mbox{for LUE},\\
      &2^{-(\alpha+\beta)/2} &&\mbox{for JUE},
    \end{aligned}
  \right. 
\end{equation}
and
\begin{equation}
  D_2(\infty) = \exp{\left(\frac{1}{2 \pi} \int\limits_{-1}^1 \frac{\tilde{f}(s)}{\sqrt{1-s^2}} \, ds \right)}.
\end{equation}

Now, we change variables in~\ref{problem_n-rh}:
\begin{equation}
  \label{glob_1}
  C(z) = (D(\infty))^{-\sigma_3}N(z)(D(z))^{\sigma_3}
\end{equation}
and find that~$C(z)$ satisfies the following RH problem.
\begin{problem}[C-RH]
  \namedlabel{problem_c-rh}{Problem~C-RH}
  \leavevmode
  \begin{enumerate}[label=\textnormal{({\roman*})},ref=C-RH-{\roman*}]
  \item $C(z)$ is analytic in~$\mathbb{C}\setminus[-1,1]$;
  \item $C^+(x) = C^-(x)   \begin{bmatrix}
      0 & 1\\
      -1 & 0
    \end{bmatrix},\ x \in (-1,1)$;
  \item $C(z) =  (I + \mathcal{O}(1/z))$ as~$z \to \infty$.
  \end{enumerate}
\end{problem}

The jump matrix in~\ref{problem_c-rh} is constant, and thus it is easy to solve this problem by diagonalizing the matrix. The solution reads
\begin{equation}
  \label{C_matr}
  \begin{aligned}
    C(z) &=    
    \begin{bmatrix}
      1 & i\\
      i & 1
    \end{bmatrix}
    (q(z))^{\sigma_3}
    \begin{bmatrix}
      1 & i\\
      i & 1
    \end{bmatrix}^{-1}\\
    &=
    \begin{bmatrix}
      \frac{1}{2}(q(z) + q^{-1}(z)) & \frac{1}{2 i}(q(z) - q^{-1}(z))\\
      -\frac{1}{2i}(q(z) - q^{-1}(z)) & \frac{1}{2}(q(z) + q^{-1}(z))
    \end{bmatrix},
  \end{aligned}
\end{equation}
where~$q(z) = (\frac{z-1}{z+1})^{1/4}$ and the principal branch of the root is used.

Going back to the original variable~$N(z)$ by using~\eqref{glob_1}, we obtain the global parametrix,
\begin{equation}
  \label{global_parametrix}
  \begin{aligned}
    N(z) &= (D(\infty))^{\sigma_3} C(z) (D(z))^{-\sigma_3}\\
    &=
    \begin{bmatrix}
      \frac{ D(\infty)}{2 D(z)}(q(z) + q^{-1}(z)) & \frac{D(z) D(\infty)}{2 i}(q(z) - q^{-1}(z))\\
      -\frac{1}{2i D(z)  D(\infty)}(q(z) - q^{-1}(z)) & \frac{D(z)}{2 D(\infty)}(q(z) + q^{-1}(z))
    \end{bmatrix}.
  \end{aligned}
\end{equation}

\begin{remark}
  \label{remark:nonuniq_glob}
  It is worth noticing that due to the fact that the asymptotic behavior at~$z=-1$ and~$z=1$ is not specified, the multiplicative RH problem~\eqref{RH_mult_prob} has several solutions besides~\eqref{szego_func}.
\end{remark}
At the end, we indicate that~$\det{N(z)}=1$, and, in particular, the matrix~$N(z)$ is invertible. Also, taking the derivative, we can easily check
\begin{equation}
  \label{Ninv_Nprime}
  N^{-1}(z)N'(z) =
  \begin{bmatrix}
    -\frac{D'(z)}{D(z)} & \frac{D^2(z) q'(z)}{i q(z)}\\
    -\frac{D^{-2}(z)q'(z)}{i q(z)} & \frac{D'(z)}{D(z)}
  \end{bmatrix}=
  \begin{bmatrix}
    -\frac{D'(z)}{D(z)} & \frac{D^2(z)}{2 i (z^2-1)}\\
    -\frac{D^{-2}(z)}{2i (z^2-1)} & \frac{D'(z)}{D(z)}
  \end{bmatrix},
\end{equation}
which we need for future reference.

\subsection{Local parametrix at~$z=1$}
\label{subsect_locp1}
\subsubsection{The case of GUE and LUE: a soft edge}
\label{subsect_locp1_soft}
In this section, we construct the local parametrix, a solution of~\ref{problem_t-rh} in a small neighborhood~$\Omega_1$ of~$z = 1$. Assume that we deal with GUE or LUE; consequently, $z=1$ is a soft edge.
Consider the following RH problem.
\begin{problem}[$\mathbf{P_{1,s}}$-RH]
  \namedlabel{problem_p1-rh}{Problem~$\mathrm{P_{1,s}}$-RH}
  \leavevmode
  \begin{enumerate}[label=\textnormal{({\roman*})},ref=$\mathrm{P}_{1,s}$-RH-{\roman*}]
  \item \label{p1rh_cond1} $P_{1,s}(z)$ is analytic in~$\Omega_1 \setminus L$;
  \item \label{p1rh_cond2} $P_{1,s}^+(z) = P_{1,s}^-(z) J_T(z),\, z \in L\cap \Omega_1$, where~$J_T(z)$ is defined in~\eqref{Y2RH};
  \item \label{p1rh_cond3} $P_{1,s}(z) (N(z))^{-1}= I + \mathcal{O}(1/n)$ as~$n \to \infty$, uniformly on~$\partial{\Omega}_1$;
  \item \label{p1rh_cond4} $P_{1,s}(z)$ is bounded at~$z = 1$.
  \end{enumerate}
\end{problem}
The conditions~\eqref{p1rh_cond1}, \eqref{p1rh_cond2}, and~\eqref{p1rh_cond4} ensure that, locally, the function~$P_{1,s}(z)$ behaves just like~$T(z)$. The condition~\eqref{p1rh_cond3} makes sure that~$P_{1,s}(z)$ matches with the global parametrix $N(z)$ on the boundary~$\partial{\Omega}_1$ up to the term of order~$O(1/n)$ as~$n \to \infty$.

An advantage of~\ref{problem_p1-rh} is that it is local and thus can be transformed into a problem with a piecewise-constant jump matrix without breaking the normalization at the infinity. Indeed, set
\begin{equation}
  \label{loc_par_p1_p1tilde}
  \tilde{P}_{1,s}(z) = P_{1,s}(z) e^{-\frac{n}{2} \phi(z)\sigma_3} (\chi(z))^{\sigma_3/2},\quad z\in \Omega_1 \setminus L.
\end{equation}
Since~$\Omega_1$ is arbitrarily small, we can always think that~$\chi(z)$ is analytic in~$\Omega_1 \setminus L$. It is easy to check that~$\tilde{P}_{1,s}(z)$ satisfies the following RH problem.
\begin{problem}[$\mathbf{\tilde{P}_{1,s}}$-RH]
  \namedlabel{problem_p1tld-rh}{Problem~$\tilde{\mathrm{P}}_{1,s}$-RH}
  \leavevmode
  \begin{enumerate}[label=\textnormal{({\roman*})},ref=$\tilde{\mathrm{P}}_{1,s}$-RH-{\roman*}]
  \item \label{p1tldrh_cond1} $\tilde{P}_{1,s}(z)$ is analytic in~$\Omega_1 \setminus L$;
  \item \label{p1tldrh_cond2} $\tilde{P}_{1,s}^+(z) = \tilde{P}_{1,s}^-(z) J_{1,s}(z),\, z \in L\cap \Omega_1$,
    \begin{equation}
      \label{P1tldsoft-RH}
      J_{1,s}(z) =\left\{
        \begin{aligned}
          &\begin{bmatrix}
            1 & 1\\
            0 & 1
          \end{bmatrix}, &&z \in (1, +\infty) \cap \Omega_1,\\
          &\begin{bmatrix}
            1 & 0\\
            1 & 1
          \end{bmatrix}, &&z \in (L_+ \cup L_-) \cap \Omega_1,\\
          &\begin{bmatrix}
            0 & 1\\
            -1 & 0
          \end{bmatrix}, &&z \in (-1,1) \cap \Omega_1;
        \end{aligned}
      \right.
    \end{equation}
  \item \label{p1tldrh_cond3} $\tilde{P}_{1,s}(z) (\chi(z))^{-\sigma_3/2} e^{\frac{n}{2} \phi(z)\sigma_3}  (N(z))^{-1}= I + \mathcal{O}(1/n)$, uniformly on~$\partial{\Omega}_1$;
  \item \label{p1tldrh_cond4} $\tilde{P}_{1,s}(z)$ is bounded at~$z = 1$.
  \end{enumerate}
\end{problem}

A solution to~\ref{problem_p1tld-rh} can be given in terms of Airy functions (e.g., see~\cite{Deift1999}). Consider the contour in Fig.~\ref{Airy_cont}, and define the jump matrix (cf.~\eqref{P1tldsoft-RH})
\begin{equation}
  \label{Airy_jumps}
  J_A(\zeta) =\left\{
    \begin{aligned}
      &\begin{bmatrix}
        1 & 1\\
        0 & 1
      \end{bmatrix}, &&\zeta \in (0, +\infty),\\
      &\begin{bmatrix}
        1 & 0\\
        1 & 1
      \end{bmatrix}, &&\zeta \in (\tilde{L}_+ \cup \tilde{L}_-),\\
      &\begin{bmatrix}
        0 & 1\\
        -1 & 0
      \end{bmatrix}, &&\zeta \in (-\infty,0).
    \end{aligned}
  \right.
\end{equation}
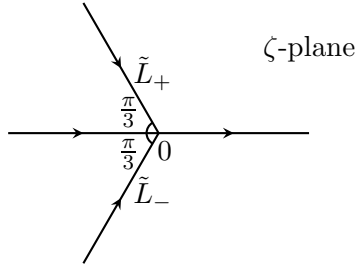
\begin{figure}[h!]
  \centering
  \begin{tikzpicture}[thick, scale=0.8]
    \begin{scope}[decoration={markings,mark= at position 0.5 with {\arrow{stealth}}}]
      \draw[postaction={decorate}] (-2.5,0) -- (0,0);
      \draw[postaction={decorate}] (0,0) -- (2.5,0);
      \draw[postaction={decorate}] (-1.25,2.165) -- (0,0);
      \draw[postaction={decorate}] (-1.25,-2.165) -- (0,0);
      
      \draw (-0.2,0) arc (180:120:0.2);
      \draw (-0.2,0) arc (-180:-120:0.2);
      \node at (2.5,1.4) {$\zeta$-plane};
      \node at (0.1,-0.3) {$0$};
      \node at (-0.5,0.35) {$\frac{\pi}{3}$};
      \node at (-0.5,-0.35) {$\frac{\pi}{3}$};
      \node at (-0.1,1) {$\tilde{L}_+$};
      \node at (-0.1,-1) {$\tilde{L}_-$};
    \end{scope}
  \end{tikzpicture}
  \caption{The auxiliary contour~$\tilde{L}=\tilde{L}_+ \cup \tilde{L}_- \cup (-\infty, +\infty)$ for the Airy parametrix.}
  \label{Airy_cont}
\end{figure}

After that, set
\begin{equation}
  \label{eq_airy_parametrix}
  A(\zeta) = \sqrt{2 \pi}
  \begin{bmatrix}
    \Ai{\zeta} & -\omega^2\Ai{\omega^2 \zeta}\\
    -i\AiPrime{\zeta} & i\omega\AiPrime{\omega^2 \zeta}
  \end{bmatrix}, \quad \arg{\zeta} \in (0, 2\pi/3),
\end{equation}
where~$\omega=e^{2\pi i /3}$,  and extend this definition to the other sectors in Fig.~\ref{Airy_cont} by applying the appropriate jumps from~\eqref{Airy_jumps}. The definition makes sense because the product of the jump matrices (taking into account the orientation of the contour) when we go around the zero is the identity matrix.

The matrix-valued function~$A(\zeta)$ constructed in this way is analytic in~$\mathbb{C} \setminus \tilde{L}$, bounded at~$\zeta = 0$, and, by definition, satisfies~$A^+(\zeta) = A^-(\zeta) J_A(\zeta),\, \zeta \in \mathbb{C} \setminus \tilde{L}$. Also, $\det{A(\zeta)}=1$, and in particular, the matrix~$A(\zeta)$ is invertible.
\begin{remark}
  \label{remark:non-uniq_Airy}
  We stress that the function satisfying the above conditions is by no means unique since its behavior as~$\zeta \to \infty$ is not specified.
\end{remark}

Using the large-$\zeta$ asymptotics of the Airy function, one can show that
\begin{equation}
  \label{Airy_Asympt}
  A(\zeta) = \zeta^{-\sigma_3/4}\frac{1}{\sqrt{2}}
  \begin{bmatrix}
    1 & i\\
    i & 1
  \end{bmatrix}
  \left(I + \mathcal{O}(\zeta^{-3/2})\right)e^{-\frac{2}{3}\zeta^{3/2}\sigma_3}
\end{equation}
as~$\zeta \to \infty$ for~$\zeta \in \mathbb{C} \setminus \tilde{L}$, where the principal branches of the roots are used.

Next, we find a conformal map~$\xi_n(z)$ from~$\Omega_1$ onto a neighborhood of~$\zeta=0$ and seek~$\tilde{\mathrm{P}}_{1,s}(z)$ in the form
\begin{equation}
  \label{P1AUXsol}
  \tilde{P}_{1,s}(z) = E_n(z)A(\xi_n(z)),
\end{equation}
where~$E_n(z)$ is some analytic function in~$\Omega_1$. Clearly, for any conformal map~$\xi_n(z)$, we have that~\eqref{P1AUXsol} satisfies~\eqref{p1tldrh_cond1}, \eqref{p1tldrh_cond2}, and \eqref{p1tldrh_cond4}; we choose this map and~$E_n(z)$ so that~\eqref{p1tldrh_cond3} is also satisfied.

Let us define~$\xi_n(z)$ in such a way that the exponent in~\eqref{p1tldrh_cond3} matches that in~\eqref{Airy_Asympt},
\begin{equation}
  \frac{2}{3}(\xi_n(z))^{3/2} = \frac{n}{2}\phi(z).
\end{equation}
The solution of this equation is
\begin{equation}
  \xi_n(z) =(3n/4)^{2/3}(\phi(z))^{2/3},
\end{equation}
where the right-hand side is analytically continued to~$\Omega_1$ by using the principal branch of the power function. Expanding~\eqref{phi} in the series about~$z=1$, we immediately see that
\begin{equation}
  \label{conf_map}
  \xi_n(z) =n^{2/3}(z-1) G(z),
\end{equation}
where~$G(z)$ is analytic in~$\Omega_1$ and~$G(1)\ne 0$. It is worth noticing that the asymptotic behavior of~$\phi(z)$ near~$z=1$ is the same for GUE and LUE; this is why we can analyze both ensembles simultaneously.

The identity~\eqref{conf_map} shows that~$\xi_n(z)$ is indeed a conformal map from~$\Omega_1$ onto some neighborhood~$\tilde{\Omega}$ of~$\zeta=0$. We note that~$(1,+\infty) \cap \Omega_1$ and~$(-1,1) \cap \Omega_1$ (with the orientation specified in Fig.~\ref{T_cont}) are mapped into~$(0,+\infty) \cap \tilde{\Omega}$ and~$(-\infty,0) \cap \tilde{\Omega}$ (with the orientation specified in Fig.~\ref{Airy_cont}). The freedom to deform the lips~$L_\pm$ can be used to ensure that the lips are mapped into~$\tilde{L}_\pm$.

Now, we choose~$\xi_n(z)$ and~$E_n(z)$ so that~\eqref{p1tldrh_cond3} is satisfied. If~$z$ is fixed and $n \to \infty$, we have~$\xi_n(z) \to \infty$, and thus the behavior of~$A(\zeta)$ as~$\zeta \to \infty$ is of relevance. Further, using~\eqref{loc_par_p1_p1tilde}, \eqref{Airy_Asympt}, \eqref{P1AUXsol}, and~\eqref{conf_map}, we see that
\begin{equation}
  \label{par1_match_cond}
  \begin{aligned}    
    P_{1,s}(z) (N(z))^{-1}=&\frac{E_n(z)}{\sqrt{2}}(3n\phi(z)/4)^{-\sigma_3/6}
    \begin{bmatrix}
      1 & i\\
      i & 1
    \end{bmatrix}\\
    &\times \bigg(I + \mathcal{O}\big( (\xi_n(z))^{-3/2} \big)\bigg) (\chi(z))^{-\sigma_3/2} (N(z))^{-1}.    
  \end{aligned}
\end{equation}
By comparing with~\eqref{p1rh_cond3}, we set
\begin{equation}
  \label{prefactor_1}
  \begin{aligned}
    E_n(z) = &N(z) (\chi(z))^{\sigma_3/2}
    \begin{bmatrix}
      \frac{1}{\sqrt{2}} & \frac{i}{\sqrt{2}}\\
      \frac{i}{\sqrt{2}} & \frac{1}{\sqrt{2}}
    \end{bmatrix}^{-1} (3n\phi(z)/4)^{\sigma_3/6}.
  \end{aligned}
\end{equation}
And the formula~\eqref{par1_match_cond} becomes 
\begin{equation}
  \label{par1_match_cond_fin}
  \begin{aligned}
    P_{1,s}(z) (N(z))^{-1}=I + &N(z) (\chi(z))^{\sigma_3/2}\mathcal{O}\big( (\xi_n(z))^{-3/2} \big)\\
    &\times  (\chi(z))^{-\sigma_3/2} (N(z))^{-1}.
  \end{aligned}
\end{equation}
Moreover, as long as~$N(z)$ and~$\chi(z)$ are uniformly bounded and
\begin{equation}
  \label{parametrics1_error}
  \mathcal{O}\big( (\xi_n(z))^{-3/2} \big) = \mathcal{O}(1/n)
\end{equation}
as~$n \to \infty$, uniformly in~$z \in \partial{\Omega}_1$, the condition~\eqref{p1rh_cond3} is satisfied.

It remains to check that~$E_n(z)$ is analytic in~$\Omega_1$. Clearly, $E_n(z)$ is analytic in~$\Omega_1 \setminus (-1,1)$, therefore it suffices to verify that there is no jump over~$(-1,1) \cap \Omega_1$ and that there is no singularity at~$z=1$. First, we use~\eqref{nrh_cond2} and compare the limits~$E_n^+(z)$ and~$E_n^-(z)$ from the above and below of~$(-1,1) \cap \Omega_1$:
\begin{equation}
  \begin{aligned}
    E_n^+(z) =& N^+(z) (\chi(z))^{\sigma_3/2}
    \begin{bmatrix}
      \frac{1}{\sqrt{2}} & \frac{i}{\sqrt{2}}\\
      \frac{i}{\sqrt{2}} & \frac{1}{\sqrt{2}}
    \end{bmatrix}^{-1} \left[(3n\phi(z)/4)^{\sigma_3/6}\right]^+\\
    =& N^-(z)
    \begin{bmatrix}
      0 & \chi(z)\\
      -\frac{1}{\chi(z)} & 0
    \end{bmatrix} (\chi(z))^{\sigma_3/2}
    \begin{bmatrix}
      \frac{1}{\sqrt{2}} & \frac{i}{\sqrt{2}}\\
      \frac{i}{\sqrt{2}} & \frac{1}{\sqrt{2}}
    \end{bmatrix}^{-1}\\
    &\times e^{\frac{\pi i \sigma_3}{2}} \left[(3n\phi(z)/4)^{\sigma_3/6}\right]^-\\
    =&N^-(z) (\chi(z))^{\sigma_3/2}
    \begin{bmatrix}
      \frac{1}{\sqrt{2}} & \frac{i}{\sqrt{2}}\\
      \frac{i}{\sqrt{2}} & \frac{1}{\sqrt{2}}
    \end{bmatrix}^{-1} \left[(3n\phi(z)/4)^{\sigma_3/6}\right]^- \\
    =& E_n^-(z).
  \end{aligned}
\end{equation}
Consequently, $E_n(z)$ is indeed analytic in~$\Omega_1\setminus \{ 1\}$ and thus can only have an isolated singularity at~$z=1$. From the explicit formula~\eqref{prefactor_1}, the order of this singularity is at most~$1/2$, and hence the singularity has to be removable. So, $E_n(z)$ is analytic in~$\Omega_1$. Also, we note that since~$\det{N(z)}=1$, one has~$\det{E_n(z)}=1$, and thus~$\det{\tilde{P}_{1,s}(z)} = \det{P_{1,s}(z)}=1$. In particular, all these matrices are non-singular.

\subsubsection{The case of JUE: a hard edge}
\label{subsubsect_bessel_param_1hard}
It remains to treat the case of JUE, where~$z=1$ is a hard edge. Again, let~$\Omega_1$ be a small neighborhood of~$z=1$. The construction of the parametrix is very similar to that in the case of the soft edge. The difference is, however, that the contour~$L \cap \Omega_1$ is different (see Fig.~\ref{T_cont}) and that the parametrix can be unbounded at~$z=1$, as it follows from~\eqref{yrh_cond4}.

Let~$P_{1,h}(z)$ be the solution of~\ref{problem_p1-rh} where instead of the condition~\eqref{p1rh_cond4} we use~\eqref{trh_cond4} for JUE. Similarly to~\eqref{Omega_omega}, define
\begin{equation}
  \tilde{\omega}(z) = \alpha \log{(z+1)} + \beta \log{(z-1)},\quad z \in \mathbb{C} \setminus (-\infty,1],
\end{equation}
where the principal branch of the logarithms is used. Then, define (cf.~\eqref{chi})
\begin{equation}
  \tilde{\chi}(z) = e^{\tilde{f}(z) + \tilde{\omega}(z)},
\end{equation}
which is clearly analytic in~$\Omega_1 \setminus (-1,1)$, and change the variables
\begin{equation}
  \label{ch_var2_hard}
  \tilde{P}_{1,h}(z) = P_{1,h}(z) e^{-\frac{n}{2} \phi(z)\sigma_3} (\tilde{\chi}(z))^{\sigma_3/2},\quad z\in \Omega_1 \setminus L.
\end{equation}
It is readily verified that~$\tilde{P}_{1,h}(z)$ solves the following RH problem with the piecewise-constant jump matrix.
\begin{problem}[$\mathbf{\tilde{P}_{1,h}}$-RH]
  \namedlabel{problem_p1htld-rh}{Problem~$\tilde{\mathrm{P}}_{1,h}$-RH}
  \leavevmode
  \begin{enumerate}[label=\textnormal{({\roman*})},ref=$\tilde{\mathrm{P}}_{1,h}$-RH-{\roman*}]
  \item \label{p1tldh_cond1} $\tilde{P}_{1,h}(z)$ is analytic in~$\Omega_1 \setminus L$;
  \item \label{p1tldh_cond2} $\tilde{P}_{1,h}^+(z) = \tilde{P}_{1,h}^-(z) J_{1,h}(z),\, z \in L\cap \Omega_1$,
    \begin{equation}
      \label{P1tldhard-RH}
      J_{1,h}(z) =\left\{
        \begin{aligned}
          &\begin{bmatrix}
            1 & 0\\
            e^{\pi i \beta} & 1
          \end{bmatrix}, &&z \in L_+ \cap \Omega_1,\\
          &\begin{bmatrix}
            1 & 0\\
            e^{-\pi i \beta} & 1
          \end{bmatrix}, &&z \in L_- \cap \Omega_1,\\
          &\begin{bmatrix}
            0 & 1\\
            -1 & 0
          \end{bmatrix}, &&z \in (-1,1) \cap \Omega_1;
        \end{aligned}
      \right.
    \end{equation}
  \item \label{p1tldh_cond3} $\tilde{P}_{1,h}(z) (\tilde{\chi}(z))^{-\sigma_3/2} e^{\frac{n}{2} \phi(z)\sigma_3}  (N(z))^{-1}= I + \mathcal{O}(1/n)$ as~$n \to \infty$, uniformly on~$\partial{\Omega}_1$;
  \item \label{p1tldh_cond4} The behavior of~$\tilde{P}_{1,h}(z)$ as~$z \to 1$:
    \begin{equation}
      \tilde{P}_{1,h}(z) =\left\{
        \begin{aligned}
          &\mathcal{O}(1) |z-1|^{\beta \sigma_3/2}, &&\beta > 0,\\
          &\mathcal{O}(|z-1|^{\beta/2}), && \beta<0,\\
          &\begin{bmatrix}
            O(1) & O(\log{|z-1|})\\
            O(1) & O(\log{|z-1|})
          \end{bmatrix}, &&\beta =0.
        \end{aligned}
      \right.
    \end{equation}
  \end{enumerate}
\end{problem}

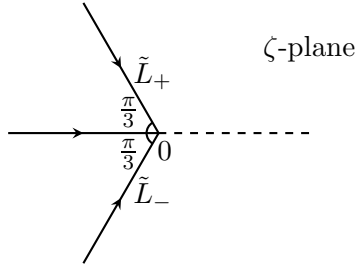
\begin{figure}[h!]
  \centering
  \begin{tikzpicture}[thick, scale=0.8]
    \begin{scope}[decoration={markings,mark= at position 0.5 with {\arrow{stealth}}}]
      \draw[postaction={decorate}] (-2.5,0) -- (0,0);
      \draw[dashed] (2.5,0) -- (0,0);
      \draw[postaction={decorate}] (-1.25,2.165) -- (0,0);
      \draw[postaction={decorate}] (-1.25,-2.165) -- (0,0);
      
      \draw (-0.2,0) arc (180:120:0.2);
      \draw (-0.2,0) arc (-180:-120:0.2);
      \node at (2.5,1.4) {$\zeta$-plane};
      \node at (0.1,-0.3) {$0$};
      \node at (-0.5,0.35) {$\frac{\pi}{3}$};
      \node at (-0.5,-0.35) {$\frac{\pi}{3}$};
      \node at (-0.1,1) {$\tilde{L}_+$};
      \node at (-0.1,-1) {$\tilde{L}_-$};

    \end{scope}
  \end{tikzpicture}
  \caption{The auxiliary contour~$\tilde{L} = \tilde{L}_+ \cup \tilde{L}_- \cup (-\infty,0)$ for the Bessel parametrix.}
  \label{Bessel_cont_hard}
\end{figure}

A solution to this problem can be given in terms of the modified Bessel functions~$I_\beta$ and~$K_\beta$ (see~\cite{Charlier2019, Vanlessen2005}). Just like in the previous section, we start by constructing a piecewise-analytic function~$\Psi_\beta(\zeta)$ such that its jump (cf.~\eqref{P1tldhard-RH}) is given by 
\begin{equation}
  \label{Psi-jump}
  J_{\Psi_\beta}(\zeta) =\left\{
    \begin{aligned}
      &\begin{bmatrix}
        1 & 0\\
        e^{\pi i \beta} & 1
      \end{bmatrix}, &&z \in \tilde{L}_+,\\
      &\begin{bmatrix}
        1 & 0\\
        e^{-\pi i \beta} & 1
      \end{bmatrix}, &&z \in \tilde{L}_-,\\
      &\begin{bmatrix}
        0 & 1\\
        -1 & 0
      \end{bmatrix}, &&z \in (-\infty,0)
    \end{aligned}
  \right.
\end{equation}
on the contour in Fig.~\ref{Bessel_cont_hard}.

Set~$\Psi_{\beta}(\zeta)$ to be
\begin{equation}
  \label{Bessel_sol_hard}
  \Psi_{\beta}(\zeta) =
  \begin{bmatrix}
    I_{\beta}(2\sqrt{\zeta}) & \frac{i}{\pi}K_{\beta}(2\sqrt{\zeta})\\
    2 \pi i \sqrt{\zeta} I_{\beta}'(2\sqrt{\zeta}) & -2 \sqrt{\zeta} K_{\beta}'(2 \sqrt{\zeta})
  \end{bmatrix},
\end{equation}
for~$\arg{\zeta} \in (-2\pi/3, 2\pi/3)$ and extend this definition to the other sectors in Fig.~\ref{Bessel_cont_hard} by applying the appropriate jumps from~\eqref{P1tldhard-RH}. It can be verified that this definition is correct, that is, the jump matrix over~$(-\infty,0)$ is consistent with the behavior of the modified Bessel functions (for details, see~\cite{Kuijlaars2004}).

The matrix-valued function~$\Psi_{\beta}(\zeta)$ constructed in this way is analytic in~$\mathbb{C} \setminus \tilde{L}$, satisfies the condition~\eqref{p1tldh_cond4}, and, by definition, satisfies~$\Psi_{\beta}^+(\zeta) = \Psi_{\beta}^-(\zeta) J_{\Psi_\beta}(\zeta),\,  \zeta \in \mathbb{C} \setminus \tilde{L}$. Moreover, $\det{\Psi_{\beta}(\zeta)}=1$, and in particular, the matrix~$\Psi_{\beta}(\zeta)$ is invertible.
\begin{remark}
  \label{remark:non-uniq_Bessel}
  Yet again, we stress that the function satisfying these conditions is by no means unique since its behavior as~$\zeta \to \infty$ is not specified.
\end{remark}

The asymptotic behavior of~$\Psi_{\beta}(\zeta)$ can be recovered from the known properties of the modified Bessel functions and is given by
\begin{equation}
  \label{Bessel_Asympt_hard1}
  \Psi_{\beta}(\zeta) = (2 \pi)^{-\sigma_3/2}\zeta^{-\sigma_3/4}\frac{1}{\sqrt{2}}
  \begin{bmatrix}
    1 & i\\
    i & 1
  \end{bmatrix}
  (I + \mathcal{O}(\zeta^{-1/2}))e^{2\zeta^{1/2}\sigma_3}
\end{equation}
for~$\zeta \in \mathbb{C} \setminus \tilde{L}$ as~$\zeta \to \infty$.

As before, we seek a solution~\ref{problem_p1htld-rh} in the form
\begin{equation}
  \label{P1-form_hard}
  \tilde{P}_{1,h}(z) = E_n(z)\Psi_{\beta}(\eta_n(z)),
\end{equation}
where~$\eta_n(z)$ is a conformal map of~$\Omega_1$ onto a neighborhood of~$\zeta=0$, and~$E_n(z)$ is analytic in~$\Omega_1$. Clearly, \eqref{p1tldh_cond1} and \eqref{p1tldh_cond2} are satisfied, and, having taken  into account the known behavior of the modified Bessel functions as~$\zeta \to 0$, we can easily check that~\eqref{p1tldh_cond4} is satisfied, too. It remains to choose~$E_n(z)$ and~$\eta_n(z)$ so that~\eqref{p1tldh_cond3} holds.

First, we match the exponents in~\eqref{p1tldh_cond3} and in~\eqref{Bessel_Asympt_hard1}:
\begin{equation}
  \label{eq_phi1_hard}
  e^{2 \sqrt{\eta_n(z)}} = e^{-\frac{n}{2} \phi(z)},
\end{equation}
which gives
\begin{equation}
  \label{eta_eq_hard}
  \eta_n(z) = \frac{n^2}{16} (\phi(z))^2.
\end{equation}
Expanding~\eqref{phi} in the series about~$z=1$, we see that
\begin{equation}
  \label{eta_asympt_eq_hard}
  \eta_n(z) =  \frac{n^2}{16}(z-1) G(z)
\end{equation}
for some analytic function~$G(z)$ which satisfies~$G(1) \ne 0$. Thus,~$\eta_n(z)$ is a conformal map of~$\Omega_1$ onto a neighborhood~$\tilde{\Omega}$ of~$\zeta=0$. We observe that~$(-1,1)\cap \Omega_1$ is mapped into~$(-\infty,0) \cap \tilde{\Omega}$. Besides, due to the freedom to deform the lips~$L_\pm$, we can always think that these lips are mapped into~$\tilde{L}_\pm$.

Now, we find the analytic factor~$E_n(z)$ so that~\eqref{p1tldh_cond3} is satisfied fully. Since~$\eta_n(z) \to \infty$ as~$n \to \infty$, the asymptotics for~$\Psi_{\beta}(\zeta)$ as~$\zeta \to \infty$ is of relevance. We use \eqref{P1-form_hard}, \eqref{Bessel_Asympt_hard1}, \eqref{eta_eq_hard}, and \eqref{ch_var2_hard} to write the left-hand side of~\eqref{p1tldh_cond3}, using the original variable~$P_{1,h}(z)$:
\begin{equation}
  \label{par1_match_cond_hard}
  \begin{aligned}    
    P_{1,h}(z) (N(z))^{-1} = &\frac{E_n(z)}{\sqrt{2}}(\pi n \phi(z)/2 )^{-\sigma_3/4}
    \begin{bmatrix}
      1 & i\\
      i & 1
    \end{bmatrix}\\
    &\times \bigg(I + \mathcal{O}\big((\eta_n(z))^{-1/2}\big)\bigg) (\tilde{\chi}(z))^{-\sigma_3/2}  (N(z))^{-1}.
  \end{aligned}
\end{equation}
In a similar way to the previous section, we define~$E_n(z)$ to be
\begin{equation}
  \label{P1-prefac_hard}
  E_n(z) = N(z) (\tilde{\chi}(z))^{\sigma_3/2}\frac{1}{\sqrt{2}}
  \begin{bmatrix}
    1 & -i\\
    -i & 1
  \end{bmatrix}
  (\pi n \phi(z)/2)^{\sigma_3/2}.
\end{equation}

The formula~\eqref{par1_match_cond_hard} becomes
\begin{equation}
  \label{par1_match_cond_hard_fin}
  \begin{aligned}
    P_{1,h}(z) (N(z))^{-1} = I + &N(z) (\tilde{\chi}(z))^{\sigma_3/2} \mathcal{O}\big((\eta_n(z))^{-1/2}\big)\\
    &\times(\tilde{\chi}(z))^{-\sigma_3/2}  (N(z))^{-1}.
  \end{aligned}
\end{equation}
Using the fact that~$N(z)$ and~$\tilde{\chi}(z)$ are uniformly bounded on~$\partial{\Omega}_1$ and that
\begin{equation}
  \label{parametrics1_error_hard}
  \mathcal{O}\big((\eta_n(z))^{-1/2}\big) = \mathcal{O}(1/n) 
\end{equation}
as~$n \to \infty$, uniformly in~$z \in \partial{\Omega}_1$, we finally arrive at~$P_{1,h}(z) (N(z))^{-1} = I+\mathcal{O}(1/n)$ as~$n \to \infty$, uniformly in~$z \in \partial{\Omega}_1$.

It is left to check that~$E_n(z)$ is analytic in~$\Omega_1$. By construction, ~$E_n(z)$ is analytic in~$\Omega_1\setminus (-1,1)$. Using~\eqref{nrh_cond2}, we will verify that~$E_n(z)$ has no jumps over~$(-1,1)$ approaching from the upper and lower half planes:
\begin{equation}
  \begin{aligned}
    E_n^+(z) =& \frac{1}{\sqrt{2}} N^+(z) [\tilde{\chi}(z)^{\sigma_3/2}]^+\cdot
    \begin{bmatrix}
      1 & -i\\
      -i & 1
    \end{bmatrix}
    \big[(\pi n \phi(z)/2)^{\sigma_3/2}\big]^+\\
    =&\frac{1}{\sqrt{2}} N^-(z)
    \begin{bmatrix}
      0 & \chi(z)\\
      -\frac{1}{\chi(z)} & 0
    \end{bmatrix}
    (\chi(z) e^{ \pi i \beta})^{\sigma_3/2}
    \begin{bmatrix}
      1 & -i\\
      -i & 1
    \end{bmatrix}\\
    &\times e^{\frac{\pi i \sigma_3}{2}} \big[(\pi n \phi(z)/2)^{\sigma_3/2}\big]^-\\
    =&\frac{1}{\sqrt{2}} N^-(z) [\tilde{\chi}(z)^{\sigma_3/2}]^-\cdot
    \begin{bmatrix}
      1 & -i\\
      -i & 1
    \end{bmatrix}
    \big[(\pi n \phi(z)/2)^{\sigma_3/2}\big]^- \\
    =& E_n^-(z).
  \end{aligned}
\end{equation}

Consequently, $E_n(z)$ is analytic in~$\Omega_1 \setminus \{ 1\}$ and can only have an isolated singularity at~$z=1$. The explicit formula~\eqref{P1-prefac_hard} shows that the order of the singularity is at most~$1/2$, and thus this singularity is removable. Consequently, $E_n(z)$ is analytic in~$\Omega_1$. Moreover, we observe that~$\det{E_n(z)}=1$, and hence~$\det{\tilde{P}_{1,h}(z)} = \det{P_{1,h}(z)} = 1$. In particular, all these matrices are non-singular.

\subsection{Local parametrix at~$z=-1$}
\label{subsect_locm1}
\subsubsection{The case of GUE: a soft edge}
The treatment of the edge~$z=-1$ virtually copies that of~$z=1$, we will only focus on the important differences. Assume that we deal with GUE, and thus~$z=-1$ is a soft edge.

Let~$\Omega_{-1}$ be a small neighborhood of~$z=-1$ such that~$\chi(z)$ is analytic in~$\Omega_{-1} \setminus L$. Consider the following RH problem.
\begin{problem}[$\mathbf{P_{-1,s}}$-RH]
  \namedlabel{problem_pm1-rh}{Problem~$\mathrm{P_{-1,s}}$-RH}
  \leavevmode
  \begin{enumerate}[label=\textnormal{({\roman*})},ref=$\mathrm{P}_{-1,s}$-RH-{\roman*}]
  \item \label{pm1s_cond1} $P_{-1,s}(z)$ is analytic in~$\Omega_{-1} \setminus L$;
  \item \label{pm1s_cond2} $P_{-1,s}^+(z) = P_{-1,s}^-(z) J_T(z),\, z \in L\cap \Omega_{-1}$, where~$J_T(z)$ is defined in~\eqref{Y2RH};
  \item \label{pm1s_cond3} $P_{-1,s}(z) (N(z))^{-1}= I + \mathcal{O}(1/n)$ as~$n \to \infty$, uniformly on~$\partial{\Omega}_{-1}$;
  \item \label{pm1s_cond4} $P_{-1,s}(z)$ is bounded at~$z = -1$.
  \end{enumerate}
\end{problem}

Locally, the function~$P_{-1,s}(z)$ has the same behavior as~$T(z)$. Again, we are going to match~$P_{-1,s}(z)$ asymptotically with the global parametrix~$N(z)$ on the boundary~$\partial{\Omega}_{-1}$ as~$n \to \infty$. Set
\begin{equation}
  \label{loc_par_p-1_p1tilde}
  \tilde{P}_{-1,s}(z) = P_{-1,s}(z) e^{-\frac{n}{2} \phi(z)\sigma_3} (\chi(z))^{\sigma_3/2},\quad z\in \Omega_{-1} \setminus L.
\end{equation}
Then~$\tilde{P}_{-1,s}(z)$ satisfies the following RH problem with a piecewise-constant jump matrix.
\begin{problem}[$\mathbf{\tilde{P}_{-1,s}}$-RH]
  \namedlabel{problem_pm1tld-rh}{Problem~$\mathrm{\tilde{P}_{-1,s}}$-RH}
  \leavevmode
  \begin{enumerate}[label=\textnormal{({\roman*})},ref=$\tilde{\mathrm{P}}_{-1,s}$-RH-{\roman*}]
  \item \label{pm1tlds_cond1} $\tilde{P}_{-1,s}(z)$ is analytic in~$\Omega_{-1} \setminus L$;
  \item \label{pm1tlds_cond2} $\tilde{P}_{-1,s}^+(z) = \tilde{P}_{-1,s}^-(z) J_{-1,s}(z),\, z \in L\cap \Omega_{-1}$,
    \begin{equation}
      \label{P-1tldsoft-RH}
      J_{-1,s}(z) =\left\{
        \begin{aligned}
          &\begin{bmatrix}
            1 & 1\\
            0 & 1
          \end{bmatrix}, &&z \in (-\infty,-1) \cap \Omega_{-1},\\
          &\begin{bmatrix}
            1 & 0\\
            1 & 1
          \end{bmatrix}, &&z \in (L_+ \cup L_-) \cap \Omega_{-1},\\
          &\begin{bmatrix}
            0 & 1\\
            -1 & 0
          \end{bmatrix}, &&z \in (-1,1) \cap \Omega_{-1};
        \end{aligned}
      \right.
    \end{equation}
  \item \label{pm1tlds_cond3} $\tilde{P}_{-1,s}(z) (\chi(z))^{-\sigma_3/2} e^{\frac{n}{2} \phi(z)\sigma_3}  (N(z))^{-1}= I + \mathcal{O}(1/n)$ as~$n \to \infty$, uniformly on~$\partial{\Omega}_{-1}$;
  \item \label{pm1tlds_cond4} $\tilde{P}_{-1,s}(z)$ is bounded at~$z = -1$.
  \end{enumerate}
\end{problem}

The parametrix~$\tilde{P}_{-1,s}(z)$ is constructed in a similar way as in Section~\ref{subsect_locp1_soft}; however, since the contour (see Fig.~\ref{T_cont}) in the vicinity of~$z=1$ is different from that in the vicinity of~$z=-1$, we need to carry out an additional transformation.

First, we reverse the orientation of the contour in Fig.~\ref{Airy_cont}. It can be checked directly that
\begin{equation}
  \label{softm1_new_A}
  \tilde{A}(\zeta) := \sigma_3A(\zeta)\sigma_3
\end{equation}
has the jump matrix~\eqref{Airy_jumps} on this reversed contour.

Next, in the same way as before, we seek a solution to~\ref{problem_pm1tld-rh} in the form
\begin{equation}
  \label{P-1AUXsol}
  \tilde{P}_{-1,s}(z) = E_n(z)\tilde{A}(\xi_n(z)),
\end{equation}
where~$\xi_n(z)$ is a conformal map of~$\Omega_{-1}$ onto some neighborhood of~$\zeta=0$ and~$E_n(z)$ is an analytic function in~$\Omega_{-1}$. It is clear that~\eqref{P-1AUXsol} satisfies~\eqref{pm1tlds_cond1}, \eqref{pm1tlds_cond2}, and~\eqref{pm1tlds_cond4}. Again, we choose~$\xi_n(z)$ and~$E_n(z)$ so that~\eqref{pm1tlds_cond3} is also satisfied.

Define
\begin{equation}
  \label{phi_tilde}
  \begin{aligned}
    \tilde{\phi}(z) &= -4 \int \limits_{-1}^z \sqrt{s^2-1}\ ds\\
    &=\left\{
      \begin{aligned}
        &\phi(z)+2\pi i,\, \Im{z}>0,\\
        &\phi(z)-2\pi i,\, \Im{z}<0,\\
      \end{aligned}
    \right.
    \quad z \in \mathbb{C}\setminus [-1,+\infty),
  \end{aligned}
\end{equation}
where~$\phi(z)$ is given in~\eqref{phi} and the principal branch of the root is used. Then set
\begin{equation}
  \label{conf_map_m1_soft_def}
  \xi_n(z) =(3n/4)^{2/3}(\tilde{\phi}(z))^{2/3},
\end{equation}
where the right-hand side is analytically continued to~$\Omega_{-1}$ by using the principal branch of the power function. Expanding~\eqref{phi_tilde} in the series about~$z=-1$, we see that
\begin{equation}
  \label{conf_map_m1_soft}
  \xi_n(z) =n^{2/3}(z+1) G(z),
\end{equation}
where~$G(z)$ is analytic in~$\Omega_{-1}$ and~$G(-1)\ne 0$.

The identity~\eqref{conf_map_m1_soft} shows that~$\xi_n(z)$ is indeed a conformal map of~$\Omega_{-1}$ onto some neighborhood~$\tilde{\Omega}$ of~$\zeta=0$. Also, we observe that~$\xi_n(z)$ maps $(-\infty,-1) \cap \Omega_{-1}$ and~$(-1,1) \cap \Omega_{-1}$ into~$(0,+\infty) \cap \tilde{\Omega}$ and~$(-\infty,0) \cap \tilde{\Omega}$. Recalling that any conformal  map preserves the angles and using the freedom to deform the lips, we see that the image of~$L \cap \Omega_{-1}$ in Fig.~\ref{T_cont} looks like the contour in Fig.~\ref{Airy_cont} with the opposite orientation. This explains why we introduced~$\tilde{A}(\zeta)$ instead of~$A(\zeta)$ in the first place.

It is clear from~\eqref{conf_map_m1_soft_def} that~$\xi_n(z) \to \infty$ as~$n \to \infty$. Therefore, again, the asymptotics of~$\tilde{A}(\zeta)$ as~$\zeta \to \infty$ is of relevance. It is easy to check that
\begin{equation}
  \label{p1mrh_ident_aux}
  e^{\frac{n}{2}\phi(z)} = (-1)^n e^{2/3 (\xi_n(z))^{3/2}},
\end{equation}
which follows from~\eqref{phi_tilde} and~\eqref{conf_map_m1_soft_def}.

Now, using the formulas~\eqref{loc_par_p-1_p1tilde}, \eqref{softm1_new_A}, \eqref{Airy_Asympt}, \eqref{P-1AUXsol}, \eqref{conf_map_m1_soft_def}, and~\eqref{p1mrh_ident_aux}, we see that the left side of the corresponding matching condition becomes
\begin{equation}
  \label{par-1_match_cond}
  \begin{aligned}
    P_{-1,s}(z) (N(z))^{-1}=&\frac{(-1)^n E_n(z)}{\sqrt{2}}(3n\tilde{\phi}(z)/4)^{-\sigma_3/6}
    \begin{bmatrix}
      1 & -i\\
      -i & 1
    \end{bmatrix}\\
    &\times \bigg(I + \mathcal{O}\big( (\xi_n(z))^{-3/2} \big)\bigg) (\chi(z))^{-\sigma_3/2} (N(z))^{-1}.    
  \end{aligned}
\end{equation}
Setting
\begin{equation}
  \label{prefactor_-1}
  \begin{aligned}
    E_n(z) = &\frac{(-1)^n N(z)}{\sqrt{2}} (\chi(z))^{\sigma_3/2}
    \begin{bmatrix}
      1 & i\\
      i & 1
    \end{bmatrix} (3n \tilde{\phi}(z)/4)^{\sigma_3/6},
  \end{aligned}
\end{equation}
we see that~\eqref{par-1_match_cond} becomes 
\begin{equation}
  \label{par-1_match_cond_fin}
  \begin{aligned}
    P_{-1,s}(z) (N(z))^{-1}=I + &N(z) (\chi(z))^{\sigma_3/2}\mathcal{O}\big( (\xi_n(z))^{-3/2} \big)\\
    &\times (\chi(z))^{-\sigma_3/2} (N(z))^{-1};
  \end{aligned}
\end{equation}
and since~$N(z)$ and~$\chi(z)$ are uniformly bounded and
\begin{equation}
  \label{parametrics-1_error}
  \mathcal{O}\big( (\xi_n(z))^{-3/2} \big) = \mathcal{O}(1/n)
\end{equation}
as~$n \to \infty$, uniformly in~$z \in \partial{\Omega}_1$, the matching condition~$P_{-1,s}(z) (N(z))^{-1} = I +\mathcal{O}(1/n)$ as~$n \to \infty$ is satisfied.

The analyticity of~$E_n(z)$ in~$\Omega_{-1}$ follows from the corresponding argument in Section~\ref{subsect_locp1_soft} mutatis mutandis. Our final observation is that since~$\det{N(z)}=1$, we also have~$\det{E_n(z)}=1$, and thus~$\det{\tilde{P}_{-1,s}(z)} = \det{P_{-1,s}(z)}=1$. In particular, all these matrices are non-singular.

\subsubsection{The case of LUE and JUE: a hard edge}
It remains to give the construction of the local parametrix in the neighborhood~$\Omega_{-1}$ of~$z=-1$ for LUE and JUE, in which case~$z=-1$ is a hard edge. The construction is very similar to that in Section~\ref{subsubsect_bessel_param_1hard}.

First, let~$P_{-1,h}(z)$ be the solution of~\ref{problem_pm1-rh}, where instead of~\eqref{pm1s_cond4} we use~\eqref{trh_cond4} for LUE and JUE. For further convenience, let us set (cf.~\eqref{Omega_omega})
\begin{equation}
  \tilde{\omega}(z)= \left\{
    \begin{aligned}
      &\alpha \log{(-z-1)} &&\mbox{for LUE},\\
      &\alpha \log{(-z-1)} + \beta \log{(1-z)} &&\mbox{for JUE},
    \end{aligned}
  \right.
\end{equation}
where~$z \in \mathbb{C} \setminus [-1,+\infty)$ and the principal branch of the logarithms is used. Also, define~$\tilde{\chi}(z)$ by
\begin{equation}
  \tilde{\chi}(z) = e^{\tilde{f}(z) + \tilde{\omega}(z)}.
\end{equation}
Clearly, the function~$\tilde{\chi}(z)$ is analytic in~$\Omega_{-1} \setminus (-1,1)$, and we can change variables:
\begin{equation}
  \label{ch_var-1_hard}
  \tilde{P}_{-1,h}(z) = P_{-1,h}(z) e^{-\frac{n}{2} \phi(z)\sigma_3} (\tilde{\chi}(z))^{\sigma_3/2},\quad z\in \Omega_{-1} \setminus L.
\end{equation}
This leads us to the RH problem with a piecewise-constant jump matrix:

\begin{problem}[$\mathbf{\tilde{P}_{-1,h}}$-RH]
  \namedlabel{problem_pm1h-rh}{Problem~$\mathrm{\tilde{P}_{-1,h}}$-RH}
  \leavevmode
  \begin{enumerate}[label=\textnormal{({\roman*})},ref=$\tilde{\mathrm{P}}_{-1,h}$-RH-{\roman*}]
  \item \label{pm1h_cond1} $\tilde{P}_{-1,h}(z)$ is analytic in~$\Omega_{-1} \setminus L$;
  \item \label{pm1h_cond2} $\tilde{P}_{-1,h}^+(z) = \tilde{P}_{-1,h}^-(z) J_{-1,h}(z),\, z \in L\cap \Omega_{-1}$,
    \begin{equation}
      \label{P-1tldh-RH}
      J_{-1,h}(z) =\left\{
        \begin{aligned}
          &\begin{bmatrix}
            1 & 0\\
            e^{- \pi i \alpha} & 1
          \end{bmatrix}, &&z \in L_+ \cap \Omega_{-1},\\
          &\begin{bmatrix}
            1 & 0\\
            e^{\pi i \alpha} & 1
          \end{bmatrix}, &&z \in L_- \cap \Omega_{-1},\\
          &\begin{bmatrix}
            0 & 1\\
            -1 & 0
          \end{bmatrix}, &&z \in (-1,1) \cap \Omega_{-1};
        \end{aligned}
      \right.
    \end{equation}
  \item \label{pm1h_cond3} $\tilde{P}_{-1,h}(z) (\tilde{\chi}(z))^{-\sigma_3/2} e^{\frac{n}{2} \phi(z)\sigma_3}  (N(z))^{-1}= I + \mathcal{O}(1/n)$ as~$n \to \infty$, uniformly on~$\partial{\Omega}_{-1}$;
  \item \label{pm1h_cond4} The behavior of~$\tilde{P}_{-1,h}(z)$ as~$z \to -1$:
    \begin{equation}
      \tilde{P}_{-1,h}(z) =\left\{
        \begin{aligned}
          &\mathcal{O}(1) |z+1|^{\alpha \sigma_3/2}, &&\alpha > 0,\\
          &\mathcal{O}(|z+1|^{\alpha/2}), && \alpha<0,\\
          &\begin{bmatrix}
            O(1) & O(\log{|z+1|})\\
            O(1) & O(\log{|z+1|})
          \end{bmatrix}, &&\alpha =0.
        \end{aligned}
      \right.
    \end{equation}
  \end{enumerate}
\end{problem}

The parametrix~$\tilde{P}_{-1,h}(z)$ is constructed in a similar manner as in Section~\ref{subsubsect_bessel_param_1hard}; however, since the contour (see Fig.~\ref{T_cont}) in the vicinity of~$z=1$ is different from that in the vicinity of~$z=-1$, we, again, need to use an additional transformation.

Reverse the orientation of the contour in Fig.~\ref{Bessel_cont_hard}. It is readily verified that
\begin{equation}
  \label{eq_psi_tilde_-1_hard}
  \tilde{\Psi}_{\alpha}(\zeta) := \sigma_3\Psi_{\alpha}(\zeta)\sigma_3
\end{equation}
has the jump matrix~\eqref{Psi-jump} on this reversed contour.

Yet again, we seek a solution to~\ref{problem_pm1h-rh} in the form
\begin{equation}
  \label{P-1AUXsol_hard}
  \tilde{P}_{-1,h}(z) = E_n(z)\tilde{\Psi}(\eta_n(z)),
\end{equation}
where~$\eta_n(z)$ is a conformal map from~$\Omega_{-1}$ onto some neighborhood of~$\zeta=0$, and~$E_n(z)$ is an analytic function in~$\Omega_{-1}$.

Define
\begin{equation}
  \label{phi_tilde_hard}
  \begin{aligned}
    \tilde{\phi}(z) &=\left\{
      \begin{aligned}
        &-2 \int \limits_{-1}^z \sqrt{\frac{s-1}{s+1}}\ ds &&\mbox{for LUE},\\
        &-2 \int \limits_{-1}^z \frac{1}{\sqrt{s^2-1}}\ ds &&\mbox{for JUE},
      \end{aligned}
    \right.\\
    &= \left\{
      \begin{aligned}
        &-\phi(z)-2\pi i,\, \Im{z}>0,\\
        &-\phi(z)+2\pi i,\, \Im{z}<0,\\
      \end{aligned}
    \right.\quad z \in \mathbb{C}\setminus [-1,+\infty),
  \end{aligned}
\end{equation}
where~$\phi(z)$ is given in~\eqref{phi} and the principal branch of the root is used.

Applying the idea from the previous sections to match the exponents in~\eqref{pm1h_cond3} and in~\eqref{Bessel_Asympt_hard1}, we set~$\eta_n(z)$ to be
\begin{equation}
  \label{eta_eq_hard_-1}
  \eta_n(z) = \frac{n^2}{16} (\tilde{\phi}(z))^2.
\end{equation}

Expanding~$\tilde{\phi}(z)$ in the series about~$z=-1$, we find that
\begin{equation}
  \label{eta_asympt_eq_hard_-1}
  \eta_n(z) =  \frac{n^2}{16}(z+1) G(z),
\end{equation}
where~$G(z)$ is analytic in~$\Omega_{-1}$ and~$G(-1) \ne 0$. Consequently, $\eta_n(z)$ is a conformal map of~$\Omega_{-1}$ onto some neighborhood~$\tilde{\Omega}$ of~$\zeta=0$, and we see that~$(-1,1)\cap \Omega_{-1}$ is mapped into~$(-\infty,0) \cap \tilde{\Omega}$. Thanks to the angle-preserving property of a conformal map and the freedom to deform the lips, it follows that the resulting contour is the one in Fig.~\ref{Bessel_cont_hard} but with the opposite orientation. Also, notice that~$L_\pm$ are mapped into~$\tilde{L}_\mp$, in particular~\eqref{pm1h_cond2} is satisfied for~\eqref{P-1AUXsol_hard} (cf.~\eqref{P-1tldh-RH} to~\eqref{P1tldhard-RH}).

Now, we notice that~$\eta_n(z) \to \infty$ as~$n \to \infty$, and thus the asymptotics of~$\Psi_{\alpha}(\zeta)$ as~$\zeta \to \infty$ is of relevance. It follows from~\eqref{phi_tilde_hard} and~\eqref{eta_eq_hard_-1} that
\begin{equation}
  \label{eq_phi-1_hard}
  e^{2 \sqrt{\eta_n(z)}} = (-1)^ne^{-\frac{n}{2} \phi(z)}.
\end{equation}
Consequently, using~\eqref{ch_var-1_hard}, \eqref{Bessel_Asympt_hard1}, \eqref{eq_psi_tilde_-1_hard}, \eqref{P-1AUXsol_hard}, \eqref{eta_eq_hard_-1}, and~\eqref{eq_phi-1_hard}, we can write the matching condition:
\begin{equation}
  \label{par-1_match_cond_hard}
  \begin{aligned}
    P_{-1,h}(z) (N(z))^{-1} = &E_n(z)(\pi n \tilde{\phi}(z)/2)^{-\sigma_3/4}\frac{(-1)^n}{\sqrt{2}}
    \begin{bmatrix}
      1 & -i\\
      -i & 1
    \end{bmatrix}
    \\
    &\times (I + \mathcal{O}(\eta_n^{-1/2}(z))) (\tilde{\chi}(z))^{-\sigma_3/2} (N(z))^{-1}.
  \end{aligned}
\end{equation}

Setting
\begin{equation}
  \label{P-1-prefac-hard}
  E_n(z) = (-1)^n N(z) (\tilde{\chi}(z))^{\sigma_3/2}\frac{1}{\sqrt{2}}
  \begin{bmatrix}
    1 & i\\
    i & 1
  \end{bmatrix}
  (\pi n \tilde{\phi}(z)/2)^{\sigma_3/2},
\end{equation}
we obtain
\begin{equation}
  \begin{aligned}
    P_{-1,h}(z) (N(z))^{-1} = I + &N(z) (\tilde{\chi}(z))^{\sigma_3/2}\mathcal{O}(\eta_n^{-1/2}(z))\\
    &\times (\tilde{\chi}(z))^{-\sigma_3/2} (N(z))^{-1}.
  \end{aligned}
\end{equation}
And since~$N(z)$ and~$\tilde{\chi}(z)$ are uniformly bounded and
\begin{equation}
  \label{parametrics-1_error_hard}
  \mathcal{O}\big((\eta_n(z))^{-1/2}\big) = \mathcal{O}(1/n) 
\end{equation}
as~$n \to \infty$, uniformly in~$z \in \partial{\Omega}_{-1}$, we finally arrive at~$P_{-1,h}(z) (N(z))^{-1} = I+\mathcal{O}(1/n)$ as~$n\to \infty$, uniformly in~$z \in \partial{\Omega}_{-1}$, as desired.

An argument similar to that in~Section~\ref{subsubsect_bessel_param_1hard} shows that~$E_n(z)$ is analytic in~$\Omega_{-1}$. Also, since~$\det{N(z)}=1$, we observe that~$\det{E_n(z)}=1$, and thus~$\det{\tilde{P}_{-1,h}(z)} = \det{P_{-1,h}(z)}=1$. In particular, all of these matrices are non-singular.

\subsection{Final transformation:  a small-norm problem}
Now, we are ready to write a small-norm problem. Let~$R(z)$ be
\begin{equation}
  \label{R-final}
  R(z)=\left\{
    \begin{aligned}
      &T(z)(P_{-1}(z))^{-1}, && z \in \Omega_{-1},\\
      &T(z)(P_{1}(z))^{-1}, && z \in \Omega_1,\\
      &T(z)(N(z))^{-1}, && z \in \mathbb{C} \setminus (\overline{\Omega_{-1} \cup \Omega_1}),
    \end{aligned}
  \right.
\end{equation}
where
\begin{equation}
  P_{1}(z):= \left\{
    \begin{aligned}
      &P_{1,s}(z) &&\mbox{for GUE, LUE},\\
      &P_{1,h}(z) &&\mbox{for JUE},
    \end{aligned}
  \right.
\end{equation}
and
\begin{equation}
  P_{-1}(z):=\left\{
    \begin{aligned}
      &P_{-1,s}(z) &&\mbox{for GUE},\\
      &P_{-1,h}(z) &&\mbox{for LUE, JUE}.
    \end{aligned}
  \right.  
\end{equation}

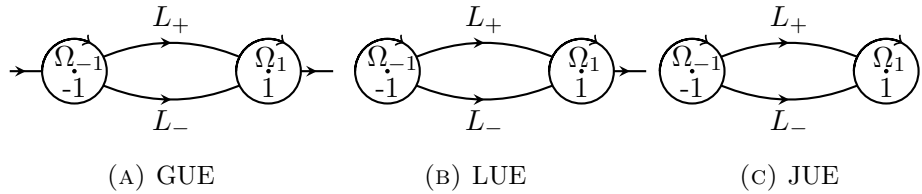
\begin{figure}[h!]
  \centering
  \begin{subfigure}[b]{0.32\textwidth}
    \centering
    \begin{tikzpicture}[thick,scale=0.86]
      \tikzset{
        clip even odd rule/.code={\pgfseteorule}, 
        invclip/.style={
          clip,insert path=
          [clip even odd rule]{
            [reset cm](-\maxdimen,-\maxdimen)rectangle(\maxdimen,\maxdimen)
          }
        }
      }
      \begin{scope}[decoration={markings,mark= at position 0.5 with {\arrow{stealth}}}]
        
        \draw [->] (-0.5,0) arc (180:-300:0.5);
        \draw [->] (2.5,0) arc (180:-300:0.5);
        
        \node at (0.08,0.19) {$\Omega_{-1}$};
        \fill (0,0) circle (1pt);
        \node at (0,-0.25) {-1};
        \node at (3.05,0.19) {$\Omega_1$};
        \fill (3,0) circle (1pt);
        \node at (3,-0.25) {1};
        \node at (1.5,0.8) {$L_+$};
        \node at (1.5,-0.8) {$L_-$};
        
        \begin{pgfinterruptboundingbox}
          \path [invclip] (0,0) circle(0.5);
          \path [invclip] (3,0) circle(0.5);
        \end{pgfinterruptboundingbox}
        
        \draw[postaction={decorate}] (3.5,0) -- (4,0);
        \draw[postaction={decorate}] (-1,0) -- (-0.5,0);
        
        \draw[postaction={decorate}] (0,0) to [out=30, in=150] (3,0);
        \draw[postaction={decorate}] (0,0) to [out=-30, in=-150] (3,0);
      \end{scope}
    \end{tikzpicture}
    \subcaption{GUE}
  \end{subfigure}
  \begin{subfigure}[b]{0.32\textwidth}
    \centering
    \begin{tikzpicture}[thick, scale =0.86]
      \tikzset{
        clip even odd rule/.code={\pgfseteorule}, 
        invclip/.style={
          clip,insert path=
          [clip even odd rule]{
            [reset cm](-\maxdimen,-\maxdimen)rectangle(\maxdimen,\maxdimen)
          }
        }
      }
      \begin{scope}[decoration={markings,mark= at position 0.5 with {\arrow{stealth}}}]
        
        \draw [->] (-0.5,0) arc (180:-300:0.5);
        \draw [->] (2.5,0) arc (180:-300:0.5);
        
        \node at (0.08,0.19) {$\Omega_{-1}$};
        \fill (0,0) circle (1pt);
        \node at (0,-0.25) {-1};
        \node at (3.05,0.19) {$\Omega_1$};
        \fill (3,0) circle (1pt);
        \node at (3,-0.25) {1};
        \node at (1.5,0.8) {$L_+$};
        \node at (1.5,-0.8) {$L_-$};
        
        \begin{pgfinterruptboundingbox}
          \path [invclip] (0,0) circle(0.5);
          \path [invclip] (3,0) circle(0.5);
        \end{pgfinterruptboundingbox}

        \draw[postaction={decorate}] (3.5,0) -- (4,0);
        \draw[draw=none] (-1,0) -- (-0.5,0);
        
        \draw[postaction={decorate}] (0,0) to [out=30, in=150] (3,0);
        \draw[postaction={decorate}] (0,0) to [out=-30, in=-150] (3,0);

      \end{scope}
    \end{tikzpicture}
    \subcaption{LUE}
  \end{subfigure}
  \begin{subfigure}[b]{0.32\textwidth}
    \centering
    \begin{tikzpicture}[thick, scale=0.86]
      \tikzset{
        clip even odd rule/.code={\pgfseteorule}, 
        invclip/.style={
          clip,insert path=
          [clip even odd rule]{
            [reset cm](-\maxdimen,-\maxdimen)rectangle(\maxdimen,\maxdimen)
          }
        }
      }
      \begin{scope}[decoration={markings,mark= at position 0.5 with {\arrow{stealth}}}]

        \draw [->] (-0.5,0) arc (180:-300:0.5);
        \draw [->] (2.5,0) arc (180:-300:0.5);

        \node at (0.08,0.19) {$\Omega_{-1}$};
        \fill (0,0) circle (1pt);
        \node at (0,-0.25) {-1};
        \node at (3.05,0.19) {$\Omega_1$};
        \fill (3,0) circle (1pt);
        \node at (3,-0.25) {1};
        \node at (1.5,0.8) {$L_+$};
        \node at (1.5,-0.8) {$L_-$};
        
        \begin{pgfinterruptboundingbox}
          \path [invclip] (0,0) circle(0.5);
          \path [invclip] (3,0) circle(0.5);
        \end{pgfinterruptboundingbox}

        \draw[postaction={decorate}] (0,0) to [out=30, in=150] (3,0);
        \draw[postaction={decorate}] (0,0) to [out=-30, in=-150] (3,0);
      \end{scope}
    \end{tikzpicture}
    \subcaption{JUE}
  \end{subfigure}
  \caption{The contour~$\Sigma_R$ of the small-norm problem.}
  \label{R_cont}
\end{figure}

By construction,~$R(z)$ only have jumps over the corresponding contour in Fig.~\ref{R_cont}.
And, as it follows from the consideration of the previous sections, this~$R(z)$ solves the following RH problem.
\begin{problem}[R-RH]
  \namedlabel{problem_r-rh}{Problem~R-RH}
  \leavevmode
  \begin{enumerate}[label=\textnormal{({\roman*})},ref=R-RH{\roman*}]
  \item $R(z)$ is analytic in~$\mathbb{C} \setminus \Sigma_R$;
  \item $R^+(z) = R^-(z) J_R(z)$, where
    \begin{equation}
      \label{jump_R}
      J_R(z) =\left\{
        \begin{aligned}
          &P_1(z) (N(z))^{-1}, &&z \in \partial \Omega_1,\\
          &P_{-1}(z) (N(z))^{-1}, &&z \in \partial \Omega_{-1},\\
          &N(z) J_T(z)(N(z))^{-1}, &&z \in L_+ \cup L_-,\\      
          &N(z) J_U(z)(N(z))^{-1}, &&z \in \mathcal{I} \setminus ([-1,1] \cup \Omega_{-1} \cup \Omega_1),
        \end{aligned}
      \right.
    \end{equation}
  \item $R(z)= I + \mathcal{O}(1/z)$ as~$z \to \infty$;
  \item $R(z)$ is bounded if~$z$ is approaching the points of self-intersection of the contour~$\Sigma_R$.
  \end{enumerate}
\end{problem}

Since~$N(z)$ is uniformly bounded, it follows that~$J_R(z) = I+\mathcal{O}(1/n)$ as~$n \to \infty$, uniformly on the contour~$\Sigma_R$. Consequently, from the theory of small-norm RH problems (e.g., see~\cite{Beresticki2017} and~\cite[Chapter~5]{Its2011}), there exist~$n_0 \in \mathbb{N}$ such that~\ref{problem_r-rh} has a unique solution for~$n>n_0$. Besides,
\begin{equation}
  \label{R-asymptotics}
  R(z) = I + \mathcal{O}(1/n),\quad R'(z) = \mathcal{O}(1/n)
\end{equation}
as~$n \to \infty$, uniformly in~$z\in\mathbb{C} \setminus \Sigma_R$. We also see, by rolling back all transformations leading to~\ref{problem_r-rh}, that~\ref{problem_y-rh} has a unique solution for~$n>n_0$.
\begin{remark}
  \label{remark:assump_a4_no_needed}
  Recall that while carrying out the steepest descent analysis of~\ref{problem_y-rh}, we did not use Assumption~\ref{assump:A4} (see Remark~\ref{remark:no_assump_A4}). In fact, it is easy to see that Assupmtion~\ref{assump:A4} is satisfied automatically for~$n>n_0$. Indeed, the existence and uniqueness of the solution~$Y(z)$ of \ref{problem_y-rh} implies that~$Y(z)$ has to match with~\eqref{eq_probly_solution}. In particular, this means that the right-hand side of~\eqref{eq_probly_solution} has to be well-defined, and thus~$H_{n-1,n}[\tilde{f}]\ne 0$ and~$H_{n,n}[\tilde{f}] \ne 0$.
\end{remark}

As a final remark, we notice that a somewhat surprising fact takes place. Even though the choice of the global and local parametrices is not unique (see Remarks~\ref{remark:nonuniq_glob}, \ref{remark:non-uniq_Airy}, \ref{remark:non-uniq_Bessel}), in the end we were able to sew them together in a consistent fashion.

\section{Proof of Lemma~\ref{lemma2_fixed_h}}
\label{lemma1_proof}
\subsection{Deformation of the weight}
\label{subsect:def_weight_proof_lemma2}
Fix~$\varepsilon>0$, and let~$h \in (-\varepsilon, \varepsilon)$. To recover the asymptotics of~$\varphi_{f,n}(h)$, as it is usually done, one may want to consider the deformation of~$ih f(x)$:
\begin{equation}
  \label{def_test_lemma2_eq1}
  \tilde{f}_t(x) = \log{\left((1-t) + t e^{i h f(x)}\right)}, \quad t \in [0,1],
\end{equation}
and the corresponding deformation of~\eqref{weight},
\begin{equation}
  \label{proof_lemma2_test_weight}
  \tilde{w}_{n,t}(x) = e^{\tilde{f}_t(x)- Q_n(x)}, \quad x \in \mathcal{I}.
\end{equation}
Then all considerations of~Section~\ref{RH_analysis} can be repeated for~$\tilde{f}:=\tilde{f}_{t}$, and the asymptotics of the solution~$Y_{n,t}(z)$ of~\ref{problem_y-rh} can be extracted. Finally, the asymptotics of~$\varphi_{f,n}(h)$ is recovered by using~\eqref{exp-via-hankel} and by integrating over~$[0,1]$ the following differential identity (e.g., see~\cite{Beresticki2017, Charlier2018, Charlier2019, Deift2014, Krasovsky2007})
\begin{equation}
  \label{diff_identity}
  \frac{\partial}{\partial t} \log{H_{n,n}[\tilde{f}_t]} = \frac{1}{2 \pi i} \int \limits_{\mathcal{I}} \left[Y^{-1}(x) Y'(x)\right]_{2,1} \frac{\partial}{\partial t} \tilde{w}(x)\, dx,
\end{equation}
where~$Y:=Y_{n,t}$ and~$\tilde{w}:= \tilde{w}_{n,t}$.

Unfortunately, there is an issue with~\eqref{def_test_lemma2_eq1}. In order for the RH approach to go through, Assumptions~\ref{assump:A1} and \ref{assump:A3} should be satisfied. In particular, $\tilde{f}_t$ should be analytic in some complex neighborhood of~$[-1,1]$ for all~$t \in [0,1]$, which is generally not the case. Thus, generally there is no corresponding neighborhood where~$J_T^\pm(z)$ (see~\eqref{matr_fact}) is analytic, which breaks the second step of the RH approach (see Section~\ref{subset_sec_transf}).

To avoid this problem, we build up on the ideas from~\cite{Deift2014}. Instead of~\eqref{def_test_lemma2_eq1}, consider the family~$\{\tilde{f}_{l,t}\}_{l=1}^q$ of deformations,
\begin{equation}
  \label{deformation1}
  \tilde{f}_{l,t}(x) = \log{\left((1-t) + t e^{\frac{i h}{q} f(x)}\right)} +\frac{ih(l-1)}{q} f(x),\quad t \in [0,1],
\end{equation}
where the principal branch of the logarithm is used; then set~$\tilde{f}:= \tilde{f}_{l,t}$. One can always choose~$q \in \mathbb{N}$ so large that
\begin{equation}
  \label{cond1}
  \left|e^{\frac{i h}{q} f(z)} -1\right| < \frac{1}{2}
\end{equation}
for~$z$ in a small complex neighborhood of~$[-1,1]$ and for all~$h \in (-\varepsilon, \varepsilon)$. Then \eqref{deformation1} is well-defined and Assumptions~\ref{assump:A1} and~\ref{assump:A3} are satisfied.

Finally, all considerations of Section~\ref{RH_analysis} carry through with~$\tilde{f}:=\tilde{f}_{l,t}$. It is also clear that all the conclusions of asymptotic nature hold uniformly in~$t \in [0,1]$.  As a result, a unique solution of~\ref{problem_r-rh}, and thus that of~\ref{problem_y-rh}, exists for~$n>n_0$, and from~\eqref{R-asymptotics} we see that
\begin{equation}
  \label{proof_lemma2_r_asympt}
  R_{n,l,t}(z) = I + \mathcal{O}(1/n),\quad R_{n,l,t}'(z) = \mathcal{O}(1/n),
\end{equation}
as~$n \to \infty$, uniformly in~$t \in [0,1]$ and~$z \in \mathbb{C} \setminus \Sigma_R$. Also, due to the explicit formula~\eqref{eq_probly_solution} (see Remark~\eqref{remark:assump_a4_no_needed}), the solution~$R_{n,l,t}(z)$, and thus~$Y_{n,l,t}(z)$, is analytic in~$t$. We equip all relevant variables with the subscripts~$l$ and~$t$ to show that these variables correspond to~$\tilde{f}:=\tilde{f}_{l,t}$. 

The next step is to integrate the differential identity~\eqref{diff_identity}. However, for this identity to hold an assumption stronger than Assumption~\ref{assump:A4} is needed. Namely, one needs the following.
\begin{assumption}
  \label{assump:A4prime}
  $H_{k,n}[\tilde{f}] \ne 0$ for all~$k=1,\dots,n$ and~$n$ large enough. 
\end{assumption}

Since~\eqref{deformation1} is complex valued, one cannot guarantee that Assumption~\ref{assump:A4prime} is satisfied. A simple argument (e.g., see~\cite{Deift2014}), however, shows that this issue can be easily avoided. Indeed, thanks to~\eqref{math-exp}, $H_{k,n}[\tilde{f}_{l,t}]$ is analytic in~$t$. Consequently, one can always choose a finite set~$\mathcal{T}_0(n)$ such that~$H_{k,n}[\tilde{f}_{l,t}] \ne 0$, for all~$t \in [0,1] \setminus \mathcal{T}_0(n)$, $k=1,\dots,n$, and~$l=1,\dots,q$. Thus, Assumption~\ref{assump:A4prime} and the differential identity hold for these~$t$. Now, we extend the differential identity to all~$t \in [0,1]$ (e.g., see~\cite{Claeys2011, Deift2014}).

Introduce the function
\begin{equation}
  \label{eq_proof_extend_t}
  S_{n,l}(t) = H_{n,n}[\tilde{f}_{l,t}] e^{-\int\limits_0^t r_{n,l}(s)\, ds}, 
\end{equation}
where~$r_{n,l}(t)$ is the right-hand side of~\eqref{diff_identity}. Clearly, from the explicit formula~\eqref{eq_probly_solution}, the expression $\left[Y^{-1}_{n,l,t}(x) Y'_{n,l,t}(x)\right]_{2,1}$ is polynomial in~$x$; therefore, $r_{n,l}(t)$ is well-defined. Since~$Y_{n,l,t}(x)$ is analytic in~$t$, by dominated convergence, the right-hand side~$r_{n,l}(t)$ is analytic, too. If the differential identity holds, then we see that~$\frac{\partial}{\partial t} S_{n,l}(t) =0$ for~$t \in [0,1] \setminus \mathcal{T}_0(n)$. And since~$S_{n,l}(t)$ is continuously differentiable in~$t$ (even analytic), it follows that~$\frac{\partial}{\partial t} S_{n,l}(t) =0$ for all~$t \in [0,1]$. We conclude that~$S_{n,l}(t)$ is constant in~$t$.

To finish the proof we need to show that~$H_{n,n}[\tilde{f}_{l,t}] \ne 0$ for all~$t \in [0,1]$, ~$l=1,\ldots,q$. First, observe that~$\tilde{f}_{1,0} = 0$, and thus, directly from~\eqref{math-exp}, $H_{n,n}[\tilde{f}_{1,0}] \ne 0$ and~$S_{n,1}(0) \ne 0$. Due to~\eqref{eq_proof_extend_t} and due to~$S_{n,l}(t)$ being constant in~$t$, we have that~$S_{n,1}(t) \ne 0$ for all~$t \in [0,1]$. Hence,~$H_{n,n}[\tilde{f}_{1,t}] \ne 0$ for~$t \in [0,1]$. In particular, $H_{n,n}[\tilde{f}_{1,1}] \ne 0$, and since $\tilde{f}_{l-1,1} = \tilde{f}_{l,0}$, we have~$H_{n,n}[\tilde{f}_{2,0}]\ne 0$. Proceeding by induction, we see that~$H_{n,n}[\tilde{f}_{l,t}] \ne 0$ for all~$t \in [0,1]$, ~$l=1,\ldots,q$. Finally, taking the derivative of~\eqref{eq_proof_extend_t} with respect to~$t$ shows that~\eqref{diff_identity} holds for all~$t\in [0,1]$ and~$l=1,\ldots,q$.

\subsection{Integration of the differential identity}
\label{subsect:integr_dif_identity}
\begin{figure}[h!]
  \centering
  \begin{tikzpicture}[thick, scale=0.8]
    \tikzset{
      clip even odd rule/.code={\pgfseteorule}, 
      invclip/.style={
        clip,insert path=
        [clip even odd rule]{
          [reset cm](-\maxdimen,-\maxdimen)rectangle(\maxdimen,\maxdimen)
        }
      }
    }
    \begin{scope}[decoration={markings,mark= at position 0.5 with {\arrow{stealth}}}]
      
      \draw[postaction={decorate},dashed] (3,0)-- (6,0);
      \draw[postaction={decorate}] (0,0) -- (3,0);
      \draw[postaction={decorate}, dashed] (-3,0) -- (0,0);
      \draw[postaction={decorate}] (0,0) to [out=50, in=130] (3,0);
      \draw[postaction={decorate}] (0,0) to [out=-50, in=-130] (3,0);
      \draw[postaction={decorate}] (-1,0) to [out=50, in=130] (4,0);
      \draw[postaction={decorate}] (-1,0) to [out=-50, in=-130] (4,0);
      
      \node at (1.6,-0.35) {$\Omega^-$};
      \node at (1.6,0.35) {$\Omega^+$};
      \node at (2.5,0.7) {$L_+$};
      \node at (2.5,-0.75) {$L_-$};
      \node at (1.5,1.3) {$\tau_+$};
      \node at (1.5,-1.5) {$\tau_-$};
      
      \fill (0,0) circle (1pt);
      \fill (4,0) circle (1pt);
      \fill (-1,0) circle (1pt);
      \node at (4.3,-0.37) {1+$\varepsilon$};
      \node at (-1.2,-0.37) {-1-$\varepsilon$};
      \node at (-0.1,-0.25) {-1};
      \fill (3,0) circle (1pt);
      \node at (3.1,-0.25) {1};
      \draw [->] ++(0:2.5) arc (180:-300:0.5);
      \draw [->] ++(0:-0.5) arc (180:-300:0.5);
    \end{scope}
  \end{tikzpicture}
  \caption{The contour to integrate the differential identity.}
  \label{int_cont}
\end{figure}
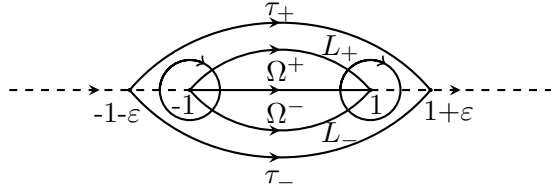

In order to integrate~\eqref{diff_identity}, we follow~\cite{Beresticki2017, Charlier2018, Charlier2019}. First, break up the contour of integration into two pieces, ~$\mathcal{I} \cap \mathcal{I}_\varepsilon$ and~$\mathcal{I} \setminus \mathcal{I}_\varepsilon$, and integrate with respect to~$t \in [0,1]$:
\begin{equation}
  \label{int_ident_1}
  \begin{aligned}
    \log{\frac{H_{n,n}\!\left[\frac{ih l}{q}f\right]}{H_{n,n}\!\left[\frac{ih (l-1)}{q}f\right]}}  = &\frac{1}{2 \pi i} \int \limits_0^1 \int \limits_{\mathcal{I} \setminus \mathcal{I}_\varepsilon} \left[Y^{-1}(x) Y'(x)\right]_{2,1} \frac{\partial}{\partial t} \tilde{w}(x) dx \, dt\\
    &+\frac{1}{2 \pi i} \int \limits_0^1 \int \limits_{\mathcal{I} \cap \mathcal{I}_\varepsilon} \left[Y^{-1}(x) Y'(x)\right]_{2,1} \frac{\partial}{\partial t} \tilde{w}(x) dx \, dt,
  \end{aligned}
\end{equation}
where~$Y:= Y_{n,l,t}$, $\tilde{w}:=\tilde{w}_{n,l,t}$, and $\mathcal{I}_{\varepsilon} = [-(1+\varepsilon), 1+\varepsilon]$. Since~$\mathcal{I} \setminus \mathcal{I}_{\varepsilon}$ is away from~$[-1,1]$ (see Fig.~\ref{int_cont}), it is possible to use the global parametrix~$N(z)$ to calculate the asymptotics of the first integral in~\eqref{int_ident_1}.

To do so, we represent~$Y^{-1}_{n,l,t}(z) Y'_{n,l,t}(z)$ via~$R_{n,l,t}(z)$ using~\eqref{U-trans}, \eqref{T-def}, and~\eqref{R-final}:
\begin{equation}
  \label{Yinv_Yprime}
  \begin{aligned}
    &Y^{-1}_{n,l,t}(z) Y'_{n,l,t}(z)\\
    &=n g'(z)\sigma_3 + e^{-n(l_R/2+g(z)) \sigma_3} N^{-1}_{l,t}(z)N'_{l,t}(z) e^{n(l_R/2+g(z)) \sigma_3}\\
    &\quad+ e^{-n(l_R/2+g(z)) \sigma_3} N^{-1}_{l,t}(z)R^{-1}_{n,l,t}(z) R'_{n,l,t}(z) N_{l,t}(z) e^{n(l_R/2+g(z)) \sigma_3}.
  \end{aligned}
\end{equation}
Note that the~$(2,1)$ matrix element of the last term is well-defined and is analytic in~$z \in \mathbb{C}\setminus [-1,1]$ because this is true for all other terms. So, we can integrate~\eqref{Yinv_Yprime} term by term.

According to~\eqref{global_parametrix}, the global parametrix~$N_{l,t}(z)$ and its inverse~$N^{-1}_{l,t}(z)$ are bounded uniformly in~$z \in \mathcal{I} \setminus \mathcal{I}_{\varepsilon}$ and $t \in [0,1]$, for all~$l=1,\dots,q$. Hence, taking into account~\eqref{proof_lemma2_r_asympt}, we see that
\begin{equation}
  \label{eq_dif_ident_est}
  N^{-1}_{l,t}(x)R^{-1}_{n,l,t}(x)R'_{n,l,t}(x)N_{l,t}(x) = \mathcal{O}(1/n)
\end{equation}
as~$n \to \infty$, uniformly in~$t \in [0,1]$ and~$x \in \mathcal{I} \setminus \mathcal{I}_{\varepsilon}$ for~$l=1,\dots,q$. Substituting~\eqref{Ninv_Nprime} into~\eqref{Yinv_Yprime}, we find that for~$x \in \mathcal{I} \setminus \mathcal{I}_{\varepsilon}$ the following holds:
\begin{equation}
  \left[Y^{-1}_{n,l,t}(x) Y'_{n,l,t}(x)\right]_{2,1} = -e^{n(l_R+2g(x)) \sigma_3} \left(\frac{D^{-2}_{l,t}(x)}{4 i x(x-1)} + O(1/n)\right), 
\end{equation}
and thus
\begin{equation}
  \label{term1}
  \begin{aligned}
    &\left[Y^{-1}_{n,l,t}(x) Y'_{n,l,t}(x)\right]_{2,1}\frac{\partial}{\partial t} \tilde{w}_{n,l,t}(x)  \\
    &= -e^{\omega(x)}\left(e^{\frac{ihl}{q}f(x)}-e^{\frac{ih(l-1)}{q}f(x)}\right) \left(\frac{D_{l,t}^{-2}(x)}{4 i x(x-1)} + O(1/n)\right)\\
    &\quad \times e^{n(l_R+2g(x)-4x)}.
  \end{aligned}
\end{equation}
Therefore, taking into account~\eqref{g_func_prop} and~\eqref{phi}, one arrives at
\begin{equation}
  \label{est_infty}
  \int \limits_{0}^{1}\int \limits_{\mathcal{I} \setminus \mathcal{I}_\varepsilon} \left[Y^{-1}_{n,l,t}(x) Y'_{n,l,t}(x)\right]_{2,1} \frac{\partial}{\partial t} \tilde{w}_{n,l,t}(x) dx\, dt=O(e^{-C n})
\end{equation}
as~$n\to \infty$, for some~$C>0$.

Now, we move on to calculating the second term in~\eqref{int_ident_1}. Observe that
\begin{equation}
  \label{eq_int_anal_transf}
  \begin{aligned}
    &\left[Y^{-1}_{n,l,t}(x) Y'_{n,l,t}(x)\right]_{2,1} \frac{\partial}{\partial t} \tilde{w}_{n,l,t}(x) \\
    &=\left(\left[Y^{-1}_{n,l,t}(x) Y'_{n,l,t}(x)\right]_{1,1}^- - \left[Y^{-1}_{n,l,t}(x) Y'_{n,l,t}(x)\right]_{1,1}^+\right)  \frac{\partial}{\partial t} \tilde{f}_{l,t}(x),
  \end{aligned}
\end{equation}
which easily follows from the direct computations using~\eqref{eq_probly_solution}. Also, it is clear that $\left[Y^{-1}_{n,l,t}(z) Y'_{n,l,t}(z)\right]_{1,1}$ is continuous over~$\mathcal{I}_{\varepsilon}\setminus \mathcal{I}$. Consequently, one has
\begin{equation}
  \label{lemma1_cont_near01}
  \begin{aligned}
    &\frac{1}{2 \pi i}\int \limits_{\mathcal{I} \cap \mathcal{I}_{\varepsilon}} \left[Y^{-1}(x) Y'(x)\right]_{2,1} \frac{\partial}{\partial t} \tilde{w}(x) \, dx \\
    &=\frac{1}{2 \pi i}\int \limits_{\mathcal{I} \cap \mathcal{I}_{\varepsilon}} \left(\left[Y^{-1}(x) Y'(x)\right]_{1,1}^- - \left[Y^{-1}(x) Y'(x)\right]_{1,1}^+\right)  \frac{\partial}{\partial t} \tilde{f}(x)\, dx \\
    &=\frac{1}{2 \pi i}\int \limits_{\mathcal{I}_{\varepsilon}} \left(\left[Y^{-1}(x) Y'(x)\right]_{1,1}^- - \left[Y^{-1}(x) Y'(x)\right]_{1,1}^+\right)  \frac{\partial}{\partial t} \tilde{f}(x)\, dx,
  \end{aligned}
\end{equation}
where~$Y:=Y_{n,l,t}$, $\tilde{w}:=\tilde{w}_{n,l,t}$, and~$\tilde{f}=\tilde{f}_{l,t}$. It is worth mentioning that for LUE and JUE, because of absence of~$\log(1\pm x)$, the quantity~$\frac{\partial}{\partial t} \tilde{f}_{l,t}(x)$ is well-defined on~$\mathcal{I}_{\varepsilon}$, while~$\frac{\partial}{\partial t} \tilde{w}_{n,l,t}(x)$ is not.

Now, by a contour deformation argument, we have
\begin{equation}
  \label{lemma1_cont_near011} 
  \begin{aligned}
    &\frac{1}{2 \pi i} \int \limits_{\mathcal{I} \cap \mathcal{I}_{\varepsilon}} \left[Y^{-1}(x) Y'(x)\right]_{2,1} \frac{\partial}{\partial t} \tilde{w}(x) \, dx \\
    &=\frac{1}{2 \pi i}\int \limits_{\mathcal{I}_{\varepsilon}} \left(\left[Y^{-1}(x) Y'(x)\right]_{1,1}^- - \left[Y^{-1}(x) Y'(x)\right]_{1,1}^+\right)  \frac{\partial}{\partial t} \tilde{f}(x)\, dx \\
    &=- \frac{1}{2 \pi i} \int \limits_{\tau_+} \left[Y^{-1}(z) Y'(z)\right]_{1,1}^+ \frac{\partial}{\partial t} \tilde{f}(z) \, dz\\
    &\quad+\frac{1}{2 \pi i}\int \limits_{\tau_-}\left[Y^{-1}(z) Y'(z)\right]_{1,1}^-  \frac{\partial}{\partial t} \tilde{f}(z)\, dz,
  \end{aligned}
\end{equation}
where~$Y:=Y_{n,l,t}$, $\tilde{w}:=\tilde{w}_{n,l,t}$, and~$\tilde{f}=\tilde{f}_{l,t}$; the contours~$\tau_\pm$ are shown in Fig.~\ref{int_cont}. Thanks to~$\tau_\pm$ being away from~$[-1,1]$,  it is again possible to use the global parametrix~$N(z)$ to find the asymptotics of the integrals.

By direct calculation in~\eqref{Yinv_Yprime} and by using~\eqref{proof_lemma2_r_asympt}, we see that
\begin{equation}
  \label{eq_diff_ident_2ndterm}
  \left[Y^{-1}_{n,l,t}(z) Y'_{n,l,t}(z)\right]_{1,1} = n g'(z) -\frac{D'_{l,t}(z)}{D_{l,t}(z)} + O(1/n),
\end{equation}
uniformly in~$t \in [0,1]$ and $z \in \tau_+ \cup \tau_-$ for~$l=1,\ldots,q$.
To deal with the first term, we use a contour deformation argument once again:
\begin{equation}
  -\frac{1}{2\pi i}\bigg(\int \limits_{\tau_+} - \int \limits_{\tau_-}\bigg) g'(z) \frac{\partial}{\partial t} \tilde{f}_{l,t}(z) dz = \int \limits_{-1}^1 \psi(x) \frac{\partial}{\partial t} \tilde{f}_{l,t}(x)\, dx,
\end{equation}
where~$\psi(x)$ is the equilibrium measure density~\eqref{eq:equil_dens}. Integrating with respect to~$t$, we find the leading term of the asymptotics
\begin{equation}
  \int \limits_0^1 \bigg[  \int \limits_{-1}^1 \psi(x) \frac{\partial}{\partial t} \tilde{f}_{l,t}(x) \, dx \bigg] \, dt = \frac{ih}{q} \int \limits_{-1}^1 f(x) \psi(x)\, dx = \frac{ih}{q} \varkappa[f].
\end{equation}

To find the next term, write
\begin{equation}
  \label{D1pD2}
  \begin{aligned}
    &\frac{1}{2\pi i}\bigg(\int \limits_{\tau_+} - \int \limits_{\tau_-}\bigg) \frac{D'_{l,t}(z)}{D_{l,t}(z)} \frac{\partial}{\partial t} \tilde{f}_{l,t}(x) \, dz\\
    &=   \frac{1}{2\pi i}\bigg(\int \limits_{\tau_+} - \int \limits_{\tau_-}\bigg) \left(\frac{D_1'(z)}{D_1(z)} + \frac{D_{2,l,t}'(z)}{D_{2,l,t}(z)}\right) \frac{\partial}{\partial t} \tilde{f}_{l,t}(x) \, dz.
  \end{aligned}
\end{equation}
It is readily verified by using~\eqref{glob_par_d1} that
\begin{equation}
  \frac{D_1'(z)}{D_1(z)} = \left\{
    \begin{aligned}
      &0 &&\mbox{for GUE},\\
      &\frac{\alpha}{2(z+1)} - \frac{\alpha}{2\sqrt{z^2-1}} &&\mbox{for LUE},\\
      &\frac{\alpha}{2(z+1)} + \frac{\beta}{2(z-1)} -\frac{\alpha+\beta}{2\sqrt{z^2-1}} &&\mbox{for JUE}.
    \end{aligned}
  \right.
\end{equation}
Then the fact that~$D_1(z)$ does not depend on~$t$ and a contour deformation argument lead to the identity
\begin{equation}
  \begin{aligned}
    &\int \limits_{0}^1 \bigg[\frac{1}{2\pi i}\bigg(\int \limits_{\tau_+} - \int \limits_{\tau_-}\bigg) \frac{D_1'(z)}{D_1(z)} \frac{\partial}{\partial t} \tilde{f}_{l,t}(z) \, dz \bigg]\, dt\\
    &= \frac{h}{2\pi q}\bigg(\int \limits_{\tau_+} - \int \limits_{\tau_-}\bigg) \frac{D_1'(z)}{D_1(z)} f(z) \, dz =\frac{i h}{q} \mu[f].
  \end{aligned}
\end{equation}

We proceed to calculate the integrals in~\eqref{D1pD2} related to~$D_{2,l,t}(z)$. In order to shorten the notation, set
\begin{equation}
  \label{theta1}
  \theta_{l,t}^{(1)}(z) = \frac{1}{2 \pi} \sqrt{z^2-1} \int\limits_{-1}^1 \frac{\log{\big((1-t) + t e^{\frac{i h}{q} f(x)}\big)}}{\sqrt{1-x^2}} \cdot \frac{dx}{z-x}
\end{equation}
and
\begin{equation}
  \label{theta2}
  \theta_{l,t}^{(2)}(z) = \frac{ih(l-1)}{2 \pi q} \sqrt{z^2-1} \int\limits_{-1}^1 \frac{f(x)}{\sqrt{1-x^2}} \cdot \frac{dx}{z-x}.
\end{equation}
Then set
\begin{equation}
  I^{(j)} =   \int \limits_0^1 \bigg[\frac{1}{2\pi i}\bigg(\int \limits_{\tau_+} - \int \limits_{\tau_-}\bigg) (\theta_{l,t}^{(j)}(z))' \frac{\partial}{\partial t} \tilde{f}_{l,t}(z) \, dz \bigg]\, dt, \quad j=1,2.
\end{equation}
Clearly, one has
\begin{equation}
  \label{int_d2}
  \int \limits_0^1 \bigg[\frac{1}{2\pi i}\bigg(\int \limits_{\tau_+} - \int \limits_{\tau_-}\bigg) \frac{D_{2,l,t}'(z)}{D_{2,l,t}(z)} \frac{\partial}{\partial t} \tilde{f}_{l,t}(z) \, dz \bigg]\, dt = I^{(1)}+I^{(2)}.
\end{equation}

To find~$I^{(1)}$, we notice that the last term in~\eqref{deformation1} does not depend on~$t$; thus, it disappears after taking the derivative~$\frac{\partial}{\partial t} \tilde{f}_{l,t}(z)$. This justifies application of Lemma~5.4 from~\cite{Beresticki2017}. We have
\begin{equation}
  I^{(1)} =  -\frac{h^2}{4 \pi^2 q^2} \int \limits_{-1}^1 \frac{f(x)}{\sqrt{1-x^2}} \mathrm{v.p.} \int \limits_{-1}^1\frac{f'(y)\sqrt{1-y^2}}{x-y} \, dy \, dx = -\frac{h^2}{2q^2} K[f].
\end{equation}

To find~$I^{(2)}$, we notice that~$\theta_{l,t}^{(2)}(z)$ does not depend on~$t$. Then applying Fubini's theorem yields
\begin{equation}
  I^{(2)} = \frac{h}{2\pi q}\bigg(\int \limits_{\tau_+} - \int \limits_{\tau_-}\bigg)  (\theta_{l,t}^{(2)}(z))' f(z) \, dz.
\end{equation}
From~\eqref{theta2}, observe that~$(\theta_{l,t}^{(2)}(z))'$ has integrable singularities at~$z=\pm1$. A contour deformation argument then yields
\begin{equation}
  \begin{aligned}
    I^{(2)} = &\frac{h}{2\pi q}\int \limits_{-1}^1\left([\theta_{l,t}^{(2)}(y)]^+ - [\theta_{l,t}^{(2)}(y)]^-\right)' f(y) \, dy.
  \end{aligned}
\end{equation}
Using integration by parts, the Sokhotski--Plemelj formulas, and properties of the Hilbert transform (e.g., see~\cite{Titchmarsh_book}), we see that
\begin{equation}
  \begin{aligned}
    I^{(2)} &= -\frac{h^2(l-1)}{2\pi^2 q^2}\int \limits_{-1}^1 \frac{f(x)}{\sqrt{1-x^2}}  \mathrm{v.p.} \int\limits_{-1}^1  \frac{f'(y) \sqrt{1-y^2}}{x-y}\, dy  \, dx \\
    &= -\frac{h^2(l-1)}{q^2} K[f].  
  \end{aligned}
\end{equation}
Finally, integrating the last term in~\eqref{eq_diff_ident_2ndterm} and collecting everything together, we arrive at the asymptotic formula
\begin{equation}
  \label{fin-assympt}
  \begin{aligned}
    \log{\frac{H_{n,n}\!\left[\frac{ih l}{q}f\right]}{H_{n,n}\!\left[\frac{ih (l-1)}{q}f\right]}} = &\frac{ih}{q} (n \varkappa[f] + \mu[f])\\
    &- \frac{h^2 (2l-1)}{2 q^2} K[f] + O\left(\frac{1}{n}\right).
  \end{aligned}
\end{equation}
Now, one can easily sum over~$l=1, \ldots, q$ and employ~\eqref{exp-via-hankel} to find that
\begin{equation}
  \label{fin_ans_lemma1}
  \begin{aligned}
    \log{\varphi_{f,n}(h)}&=\log{\mexp[n]{e^{i h\Tr{f(M)}}}}\\
    &= ih(n \varkappa[f] + \mu[f])-\frac{h^2}{2} K[f] + O\left(\frac{1}{n}\right),
  \end{aligned}
\end{equation}
which  immediately yields~\eqref{lemma2_fixed_h_form}.

The final part of the proof is to validate that~\eqref{variance_kf} can be written as~\eqref{quadr_series}. This has been done in~\cite[p.~172]{Johansson1998}.

\section{Proof of Lemma~\ref{lemma1_growing_arg}}
\label{sec_proof_lemma_gr_arg}
In order to prove the lemma, the asymptotics of~$\varphi_{f,n}(h n^\gamma)$ is needed. We follow the general procedure described in Section~\ref{lemma1_proof}; however, this time we cannot use~\eqref{deformation1}. Indeed, in our case the analog of~\eqref{deformation1} would be 
\begin{equation}
  \label{deformation1_test}
  \tilde{f}_{n,l,t}(x) = \log{\left((1-t) + t e^{\frac{i h n^\gamma}{q} f(x)}\right)} +\frac{ih n^\gamma (l-1)}{q} f(x),\quad t \in [0,1];
\end{equation}
and, since the exponent grows with~$n$, for general~$f$ there is no complex neighborhood of~$[-1,1]$ such that
\begin{equation}
  \label{deformation1_test_cond}
  (1-t) + t e^{\frac{i h n^\gamma}{q} f(z)} \ne 0
\end{equation}
for all~$n>n_0$. Thus, we cannot guarantee that Assumption~\ref{assump:A3} is satisfied.

Our idea is to delete~$n^\gamma$ from~\eqref{deformation1_test}, in this case the issue described above does not arise. On the other hand, having done so, we have to make sure that the corresponding family of deformations ends up at~$i h n^\gamma f(x)$. To achieve this, we let~$l$ be arbitrarily large and introduce the indicators~$\mathds{1}[l \le q n^{\gamma} +1]$ in our deformation.

Fix~$\varepsilon>0$, and choose~$q$ so large that
\begin{equation}
  \label{lemma1_proof_cond1_circ}
  \left|e^{\frac{i h}{q} \mathds{1}[l \le q n^{\gamma} +1] f(z)} -1 \right| < \frac{1}{2}
\end{equation}
for~$z$ in a sufficiently-small complex neighborhood of~$[-1,1]$, for all~$l, n \in \mathbb{N}$ and~$h \in (-\varepsilon, \varepsilon)$. Then introduce the deformation
\begin{equation}
  \label{deformation2}
  \begin{aligned}
    \tilde{f}_{n,l,t}(x) = &\log{\left((1-t) + t e^{\frac{i h}{q} \mathds{1}[l \le q n^{\gamma} +1] f(x)}\right)} \\
    &+\frac{ih(l-1)}{q} \mathds{1}[l \le q n^{\gamma} +1] f(x), \quad l \in \mathbb{N},
  \end{aligned}
\end{equation}
which is well-defined thanks to~\eqref{lemma1_proof_cond1_circ}.

Assumptions~\ref{assump:A1} and~\ref{assump:A3} are satisfied, and all considerations of Section~\ref{RH_analysis} carry over without change, except that, unlike~\eqref{deformation1}, the function~$\tilde{f}:=\tilde{f}_{n,l,t}$ now depends on~$n$. This means that the conclusions of asymptotic nature in the RH analysis may be affected. This means that we need to check if~\eqref{jump_R} still satisfies~$J_R(z) = I +O(1/n)$ as~$n \to \infty$ on~$\Sigma_R$ in Fig.~\ref{R_cont}.

First, let~$z \in \partial \Omega_1$. Then~\eqref{jump_R} takes on the form
\begin{equation}
  J_R(z) = P_{1;n,l,t}(z) (N_{n,l,t}(z))^{-1}.
\end{equation}
The construction of the parametrix~$P_{1;n,l,t}(z)$ requires no change; however, now~$N_{n,l,t}(z)$ depends on~$n$ and therefore may not be bounded as~$l, n \to \infty$. We need to study the right-hand side of~\eqref{par1_match_cond_fin} more carefully.

Write
\begin{equation}
  \label{par1_depl}
  \begin{aligned}
    J_R(z) = I + &N_{n,l,t}(z) (\chi_{n,l,t}(z))^{\sigma_3/2}  \mathcal{O}(1/ (n \phi(z)))\\
    &\times(\chi_{n,l,t}(z))^{-\sigma_3/2} (N_{n,l,t}(z))^{-1}, \quad z \in \partial \Omega_1,
  \end{aligned}
\end{equation}
and notice that this expression covers all three ensembles (up to a small modification for JUE, in which we need to substitute~$\chi$ with~$\tilde{\chi}$, cf.~\eqref{par1_match_cond_hard_fin}). Plugging in~\eqref{global_parametrix}, we obtain
\begin{equation}
  \label{eq_lemma1_proof_expanded_Jr}
  \begin{aligned}
    J_R(z) = I + &(D_{n,l,t}(\infty))^{\sigma_3} C(z) \left(\frac{\chi_{n,l,t}^{1/2}(z)}{D_{n,l,t}(z)}\right)^{\sigma_3}  \mathcal{O}(1/ (n \phi(z)))\\
    &\times \left(\frac{\chi_{n,l,t}^{1/2}(z)}{D_{n,l,t}(z)}\right)^{-\sigma_3}(C(z))^{-1} (D_{n,l,t}(\infty))^{-\sigma_3}, \quad z \in \partial \Omega_1.
  \end{aligned}
\end{equation}
Then using~\eqref{szego_func}, \eqref{chi}, and~\eqref{deformation2} to write out~$\frac{\chi_{n,l,t}^{1/2}(z)}{D_{n,l,t}(z)}$, we see that there is the factor
\begin{equation}
  \label{eq_our_proof_unb_shrink}
  \begin{aligned}
    &\exp{\!\Bigg(\!\!\frac{i h (l-1)}{2 q} \mathds{1}[l \le q n^\gamma + 1] \bigg(f(z) - \frac{\sqrt{z^2-1}}{\pi} \int \limits_{-1}^1 \frac{f(x) \, dx}{\sqrt{1-x^2}(z-x)}\bigg)\!\!\Bigg)}\\
    &=\exp{\Bigg( \frac{i h (l-1)}{2q} \mathds{1}[l \le q n^\gamma + 1] \frac{\sqrt{z^2-1}}{\pi} \int \limits_{-1}^1 \frac{f(z)-f(x)}{\sqrt{1-x^2}(z-x)} \, dx\Bigg)},
  \end{aligned}    
\end{equation}
which is unbounded as~$l, n \to \infty$. To compensate for the unboundedness, we notice that the integral in the right-hand side of~\eqref{eq_our_proof_unb_shrink} is uniformly bounded; therefore, making the neighborhood~$\Omega_1 := \Omega_1^{(n)}$ contract at the rate~$O(1/n^{2\gamma})$ as~$n \to \infty$ clears up the issue. So, we see that~\eqref{eq_our_proof_unb_shrink} is bounded as~$n \to \infty$, uniformly in~$z \in \partial \Omega_1^{(n)}$, $l \in \mathbb{N}$, $t\in [0,1]$, and~$h \in (-\varepsilon, \varepsilon)$.

Using the series expansion about~$z=1$ in~\eqref{phi} yields
\begin{equation}
  \label{lemma1_proof_phin_asymp}
  \mathcal{O}(1/(n \phi(z))) = \left\{
    \begin{aligned}
      &\mathcal{O}(1/n^{1-3\gamma}) &&\mbox{for GUE, LUE},\\
      &\mathcal{O}(1/n^{1-\gamma}) &&\mbox{for JUE},
    \end{aligned}
  \right. \quad
\end{equation}
as~$n \to \infty$, uniformly in~$z \in \partial \Omega_1^{(n)}$. Also, recalling the asymptotics~$C(z) = \mathcal{O}((z-1)^{-1/4})$ as~$z \to 1$ from~\eqref{C_matr} and noticing that~$D_{n,l,t}(\infty)$ from~\eqref{glob_param_D_infty} is bounded in~$l,n \in \mathbb{N}$ and~$t \in[0,1]$, we see that~\eqref{eq_lemma1_proof_expanded_Jr} turns into
\begin{equation}
  \label{P1_var_cond3}
  P_{1;n,l,t}(z) (N_{n,l,t}(z))^{-1} = I + \left\{
    \begin{aligned}
      &\mathcal{O}(1/n^{1-4\gamma}) &&\mbox{for GUE, LUE},\\
      &\mathcal{O}(1/n^{1-2\gamma}) &&\mbox{for JUE},
    \end{aligned}
  \right. \quad
\end{equation}
as~$n \to \infty$, uniformly in~$z \in \partial \Omega_1^{(n)}$, $l \in \mathbb{N}$, $t \in[0,1]$, and~$h \in (-\varepsilon, \varepsilon)$.

The local parametrix at~$z=-1$ can be handled in a similar way. We make the neighborhood~$\Omega_{-1}:=\Omega_{-1}^{(n)}$ contract at the rate~$O(1/n^{2\gamma})$ as~$n \to \infty$. Then it is readily verified that
\begin{equation}
  \label{P-1_var_cond3}
  J_R(z) = I + \left\{
    \begin{aligned}
      &\mathcal{O}(1/n^{1-4\gamma}) &&\mbox{for GUE},\\
      &\mathcal{O}(1/n^{1-2\gamma}) &&\mbox{for LUE, JUE},
    \end{aligned}
  \right. \quad
\end{equation}
as~$n \to \infty$, uniformly in~$z \in \partial \Omega_{-1}^{(n)}$, $l \in \mathbb{N}$, $t \in[0,1]$, and~$h \in (-\varepsilon, \varepsilon)$.

It remains to check if~$J_R(z)$ converges to the identity matrix on~$\mathcal{I}\setminus \left([-1,1] \cup \Omega_{-1}^{(n)} \cup \Omega_1^{(n)}\right)$ and on~$L_\pm \setminus (\Omega_{-1}^{(n)} \cup \Omega_1^{(n)})$ (see Fig.~\ref{R_cont}). From~\eqref{jump_R}, we see  that for~$x \in \mathcal{I}\setminus \left([-1,1] \cup \Omega_{-1}^{(n)} \cup \Omega_1^{(n)}\right)$ one has
\begin{equation}
  J_R(x) = N_{n,l,t}(x) J_U(x)(N_{n,l,t}(x))^{-1}.
\end{equation}
Plugging in~\eqref{jump_U} and~\eqref{global_parametrix}, we obtain
\begin{equation}
  \label{eq_lemma1_proof_jr_expr2}
  \begin{aligned}
    J_R(x) = I + &\frac{\chi_{n,l,t}(x) e^{-n \phi(x)}}{D_{n,l,t}^2(x)} (D_{n,l,t}(\infty))^{\sigma_3} C(x)
    \begin{bmatrix}
      0 & 1\\
      0 & 0
    \end{bmatrix}\\
    &\times (C(x))^{-1} (D_{n,l,t}(\infty))^{-\sigma_3}.
  \end{aligned}
\end{equation}
Using~\eqref{szego_func}, \eqref{chi}, \eqref{deformation2}, and asymptotics of~$C(x)$, we see that all factors in~\eqref{eq_lemma1_proof_jr_expr2} exhibit at most power-like growth as~$n \to \infty$, and this growth is damped by~$e^{-n \phi(x)}$. So, for some~$C>0$ we have
\begin{equation}
  \label{lemma1_proof_contr1}
  J_R(x) = I + \left\{
    \begin{aligned}
      &\mathcal{O}(e^{-C n^{1-3 \gamma}}) &&\mbox{for GUE, LUE},\\
      &\mathcal{O}(e^{-C n^{1- \gamma}}) &&\mbox{for JUE},
    \end{aligned}
  \right.
\end{equation}
as~$n \to \infty$, uniformly in~$x \in \mathcal{I}\setminus \left([-1,1] \cup \Omega_{-1}^{(n)} \cup \Omega_1^{(n)}\right)$,  $l \in \mathbb{N}$, $t \in[0,1]$, and~$h \in (-\varepsilon, \varepsilon)$.

Further, consider~$z \in L_\pm \setminus (\Omega_{-1}^{(n)} \cup \Omega_1^{(n)})$. In an analogous way as before, using~\eqref{jump_R}, \eqref{Y2RH}, and \eqref{matr_fact}, one can check that
\begin{equation}
  \label{lemma1_proof_contr2}
  \begin{aligned}
    J_R(z) &= I \mp \frac{D_{n,l,t}^2(z) e^{-n \phi(z)}}{\chi_{n,l,t}(z)} (D_{n,l,t}(\infty))^{\sigma_3} C(z)\\
      &\qquad \quad \times
      \begin{bmatrix}
        0 & 0\\
        1 & 0
      \end{bmatrix} (C(z))^{-1} (D_{n,l,t}(\infty))^{-\sigma_3}\\
    &= I + \left\{
      \begin{aligned}
        &\mathcal{O}(e^{-C n^{1-3 \gamma}}) &&\mbox{for GUE, LUE},\\
        &\mathcal{O}(e^{-C n^{1- \gamma}}) &&\mbox{for JUE},
      \end{aligned}
    \right.
  \end{aligned}
\end{equation}
uniformly in~$z \in L_\pm \setminus (\Omega_{-1}^{(n)} \cup \Omega_1^{(n)})$, $l \in \mathbb{N}$, $t \in[0,1]$, and~$h \in (-\varepsilon, \varepsilon)$.

Finally, we see that
\begin{equation}
  J_R(z) =I+ \left\{
    \begin{aligned}
      &\mathcal{O}(1/n^{1-4\gamma})&&\mbox{for GUE, LUE},\\
      &\mathcal{O}(1/n^{1-2\gamma})&&\mbox{for JUE}
    \end{aligned}
  \right.
\end{equation}
as~$n \to \infty$, on the whole contour~$\Sigma_R^{(n)}$ (same as in Fig.~\ref{R_cont} except that~$\Omega_{\pm1} := \Omega_{\pm1}^{(n)}$), uniformly in all the relevant parameters. Applying an analog of the small-norm theory but for varying (contracting) contours (e.g., see~\cite[Appendix]{Bleher2007}), we see that there exist~$n_0 \in \mathbb{N}$ such that~\ref{problem_r-rh} has a unique solution~$R_{n,l,t}(z)$ for~$n>n_0$. Besides,
\begin{equation}
  \label{R-asymptoticsl}
  R_{n,l,t}(z) =I+ \left\{
    \begin{aligned}
      &\mathcal{O}(1/n^{1-4\gamma}) &&\mbox{for GUE, LUE},\\
      &\mathcal{O}(1/n^{1-2\gamma}) &&\mbox{for JUE}
    \end{aligned}
  \right.\
\end{equation}
and
\begin{equation}
  R'_{n,l,t}(z) =\left\{
    \begin{aligned}
      &\mathcal{O}(1/n^{1-4\gamma}) &&\mbox{for GUE, LUE},\\
      &\mathcal{O}(1/n^{1-2\gamma}) &&\mbox{for JUE},
    \end{aligned}
  \right.
\end{equation}
as~$n \to \infty$, uniformly in~$l \in \mathbb{N}$, $t \in[0,1]$, $h \in (-\varepsilon, \varepsilon)$, and~$z \in \mathbb{C} \setminus \Sigma_R^{(n)}$. In a similar way as in Section~\ref{subsect:def_weight_proof_lemma2}, we see that~$R_{n,l,t}(z)$, and thus~$Y_{n,l,t}(z)$, is analytic in~$t$.
\begin{remark}
  \label{remark: loc_par_define_rate}
  It follows from~\eqref{P1_var_cond3}, \eqref{P-1_var_cond3}, \eqref{lemma1_proof_contr1}, and~\eqref{lemma1_proof_contr2} that the main contribution to the error in~\eqref{R-asymptoticsl} is due to the local parametrices. Therefore, it is the local parametrices what determine the error term in the asymptotic of~$\varphi_{f,n}(h n^\gamma)$, and thus the speed of convergence in Theorem~\ref{th1_speed_of_conv} (see the proof in Section~\ref{sect_proof_of_th1}).
\end{remark}

The next step is to integrate the differential identity. An argument similar to that of Section~\ref{subsect:def_weight_proof_lemma2} shows that the identity holds for all~$t \in [0,1]$, and we can proceed to carry out the analysis of Section~\ref{subsect:integr_dif_identity}. Note that the only conclusions affected by our choice of~$\tilde{f}_{n,l,t}$ are those of asymptotic nature: everything else carries over without serious change.

Choose the same contour of integration (independent of~$n$) as in Fig.~\ref{int_cont}. First, we calculate the integral over~$\mathcal{I} \setminus \mathcal{I}_{\varepsilon}$ in~\eqref{int_ident_1}. Since~$N_{n,l,t}(x)$ is uniformly bounded on~$\mathcal{I} \setminus \mathcal{I}_{\varepsilon}$ in~$l,n \in \mathbb{N}$, $t \in [0,1]$, and~$h \in (-\varepsilon, \varepsilon)$, we see from~\eqref{R-asymptoticsl} that (cf.~\eqref{eq_dif_ident_est})
\begin{equation}
  \label{eq_dif_ident_estl}
  N^{-1}(x)R^{-1}(x)R'(x)N(x) =\left\{
    \begin{aligned}
      &\mathcal{O}(1/n^{1-4\gamma}) &&\mbox{for GUE, LUE},\\
      &\mathcal{O}(1/n^{1-2\gamma}) &&\mbox{for JUE},
    \end{aligned}
  \right.
\end{equation}
as~$n \to \infty$, uniformly in all relevant parameters and in~$x \in \mathcal{I} \setminus \mathcal{I}_{\varepsilon}$. From here on out we often dropped the subscripts~$n, l$, and~$t$ for brevity.

By assumption, $f(x) = O(e^{A V(x)})$, $A>0$; also, for~$x \in \mathcal{I} \setminus \mathcal{I}_{\varepsilon}$ the following straightforward estimate holds
\begin{equation}
  \left|e^{\frac{i h l}{q} \mathds{1}[l \le q n^{\gamma} +1] f(x)}-e^{\frac{i h (l-1)}{q} \mathds{1}[l \le q n^{\gamma} +1] f(x)}\right| \le |h f(x)| \mathds{1}[l \le q n^{\gamma} +1].
\end{equation}
Thus, it is immediate to see that the integral in~\eqref{est_infty} becomes
\begin{equation}
  \int \limits_0^1\int \limits_{\mathcal{I} \setminus \mathcal{I}_\varepsilon} \left[Y^{-1}(x) Y'(x)\right]_{2,1} \frac{\partial}{\partial t} \tilde{w}(x) dx\, dt =\mathds{1}[l \le q n^{\gamma} +1] O(h e^{-C n})
\end{equation}
for some~$C>0$, uniformly in~$l \in \mathbb{N}$ and~$h \in (-\varepsilon,\varepsilon)$.

Next, since integrating the first two terms in~\eqref{Yinv_Yprime} carries no essential change, it remains to integrate the error term over~$\mathcal{I} \cap \mathcal{I}_{\varepsilon}$. Notice that~$N_{n,l,t}(z)$ is not bounded on~$\tau_\pm$ (see Fig.~\ref{int_cont}) as~$l,n \to \infty$ because the corresponding Szeg\H{o} function~\eqref{szego_func} is not (see~\eqref{global_parametrix}). Nevertheless, a direct calculation of the left-hand side of~\eqref{eq_dif_ident_estl} shows that $D_{n, l,t}(z)$ cancels out in the~$(1,1)$ element of this matrix. Hence, we still have
\begin{equation}
  \label{lemma1_int_last_err_term}
  \left[N^{-1}(z)R^{-1}(z)R'(z)N(z)\right]_{1,1} =\left\{
    \begin{aligned}
      &O(1/n^{1-4\gamma}) &&\mbox{for GUE, LUE},\\
      &O(1/n^{1-2\gamma}) &&\mbox{for JUE},
    \end{aligned}
  \right.
\end{equation}
as~$n \to \infty$, uniformly in all relevant parameters and in~$z \in \tau_+ \cup \tau_-$. Finally, from the straightforward estimate
\begin{equation}
  \begin{aligned}
    \left|\frac{\partial \tilde{f}(z)}{\partial t}\right| = &\left|\frac{e^{\frac{i h}{q} \mathds{1}[l \le q n^{\gamma} +1] f(z)} - 1}{(1-t) + t e^{\frac{i h}{q} \mathds{1}[l \le q n^{\gamma} +1] f(z)}}\right|\\
    &\le C|h| \mathds{1}[l \le q n^{\gamma} +1], \quad z \in \tau_\pm,
  \end{aligned}
\end{equation}
for some~$C>0$, and from~\eqref{lemma1_int_last_err_term} and~\eqref{lemma1_cont_near011}, we see that
\begin{equation}
  \begin{aligned}
    &\frac{1}{2 \pi i} \int \limits_0^1 \int \limits_{\mathcal{I} \cap \mathcal{I}_{\varepsilon}} \left[N^{-1}(x)R^{-1}(x)R'(x)N(x)\right]_{2,1} \frac{\partial}{\partial t} \tilde{w}(x) \, dx\, dt \\
    &=\mathds{1}[l \le q n^{\gamma} +1] \cdot
    \left\{
      \begin{aligned}
        &O(h/n^{1-4\gamma}) &&\mbox{for GUE, LUE},\\
        &O(h/n^{1-2\gamma}) &&\mbox{for JUE},
      \end{aligned}
    \right.
  \end{aligned}
\end{equation}
as~$n \to \infty$, uniformly in~$l \in \mathbb{N}$ and~$h \in (-\varepsilon,\varepsilon)$.

Collecting all terms, we arrive at an analogue of~\eqref{fin-assympt},
\begin{equation}
  \label{eq_148}
  \begin{aligned}
    &\log{\frac{H_{n,n}\!\left[\frac{ih l}{q}\mathds{1}[l \le q n^\gamma+1] f\right]}{H_{n,n}\!\left[\frac{ih (l-1)}{q} \mathds{1}[l \le q n^\gamma+1] f\right]}}\\
    &= \Bigg(\frac{ih}{q} (n \varkappa[f] + \mu[f]) -\frac{h^2 (2l-1)}{2 q^2} K[f] + O\left(\frac{h}{n^{1-(d-1)\gamma}}\right)\!\!\!\Bigg)\\
    &\quad\times \mathds{1}[l \le q n^\gamma+1]
  \end{aligned}
\end{equation}
as~$n \to \infty$, where~$d=5$ for GUE and LUE, $d=3$ for JUE; the~$O$-term is uniform in~$l \in \mathbb{N}$ and~$h \in (-\varepsilon, \varepsilon)$.

Summing over all~$l=1,2,\dots$, replacing~$h \mapsto h q n^\gamma/([q n^\gamma]+1)$, where~$[n^\gamma]$ is the integer part of~$n^\gamma$, and using~\eqref{exp-via-hankel} produce the final asymptotic formula
\begin{equation}
  \label{lemma_proof_form_fin_sum}
  \begin{aligned}
    \log{\varphi_{f,n}(h n^\gamma)} = \log{\frac{H_{n,n}\!\left[ih n^{\gamma} f\right]}{H_{n,n}[0]}}
    =&ih n^\gamma (n \varkappa[f] + \mu[f])\\
    &-\frac{n^{2 \gamma} h^2}{2}K[f]+ O\left(\frac{h}{n^{1-d \gamma}}\right),    
  \end{aligned}
\end{equation}
as~$n \to \infty$, uniformly in~$h \in (-\varepsilon, \varepsilon)$. If~$\gamma \le 1/d$, taking exponents of the both sides of~\eqref{lemma_proof_form_fin_sum} immediately yields the claim of the lemma.

As a final remark we note that if~$\gamma \in (1/d, 1/(d-1))$, the~$O$-term in~\eqref{lemma_proof_form_fin_sum} is growing. Nonetheless, using the inequality~$|e^z-1|\le |z| e^{|z|}$, $z \in \mathbb{C}$, one can still get
\begin{equation}
  \label{eq_big_gamma}
  \sup_{n}\sup_{|h|<\varepsilon n^\gamma} \left(\frac{n^{1-(d-1)\gamma}}{e^{C|h|}} \left| \frac{\varphi_{f,n}(h)-\varphi_{\mathcal{N}}(h)}{h\varphi_{\mathcal{N}}(h)} \right|\right) < +\infty
\end{equation}
for some~$C>0$. Notice the presence of the additional factor~$e^{C|h|}$.

\section{Proof of Theorem~\ref{th1_speed_of_conv}}
\label{sect_proof_of_th1}

We prove the theorem by means of Feller's smoothing inequality (e.g., see~\cite{Feller2,Johansson1997}). First, introduce the centered random variable
\begin{equation}
  \tilde{S}_{f,n}= S_{f,n} - \mexp[n]{S_{f,n}} = \frac{\Tr{f(M)} - \mexp[n]{\Tr{f(M)}}}{\sqrt{K[f]}},
\end{equation}
and let~$\tilde{F}_{f,n}(x)= F_{f,n}(x+\mexp[n]{S_{f,n}})$ be the corresponding cumulative distribution function. From Lemma~\ref{lemma1_growing_arg} for~$\gamma = 0$, it is easy to see that
\begin{equation}
  \mexp[n]{S_{f,n}} = O\left(\frac{1}{n}\right), \quad n \to \infty.
\end{equation}
Now, write the Kolmogorov--Smirnov distance
\begin{equation}
  \label{theorem1_proof_ineq_centering}
  \begin{aligned}
    \sup_x{|F_{f,n}(x)-F_{\mathcal{N}}(x)|} \le &\sup_x{|\tilde{F}_{f,n}(x)-F_{\mathcal{N}}(x)|} \\
    &+ \sup_x{|F_{\mathcal{N}}(x)-F_{\mathcal{N}}(x+\mexp[n]{S_{f,n}})|}
  \end{aligned}
\end{equation}
and notice that the last term is easy to estimate directly
\begin{equation}
  \sup_x{|F_{\mathcal{N}}(x)-F_{\mathcal{N}}(x+\mexp[n]{S_{f,n}})|} = 2 F_{\mathcal{N}}(|\mexp[n]{S_{f,n}}|/2)-1 = O\left(\frac{1}{n}\right)
\end{equation}
as~$n \to \infty$.

It remains to estimate~$\sup\limits_x {|\tilde{F}_{f,n}(x)-F_{\mathcal{N}}(x)|}$. Let
\begin{equation}
  \tilde{\varphi}_{f,n}(h) = \varphi_{f,n}(h) e^{-ih \mexp[n]{S_{f,n}}}
\end{equation}
be the characteristic function of~$\tilde{S}_{f,n}$. Then from Feller's inequality (see~\cite[p. 538]{Feller2}) we have the bound for~$\delta:= \sup_x{|\tilde{F}_{f,n}(x)-F_{\mathcal{N}}(x)|}$,
\begin{equation}
  \delta \le \frac{1}{\pi} \int \limits_{-T}^{T} \left|\frac{\tilde{\varphi}_{f,n}(h) - \varphi_{\mathcal{N}}(h)}{h} \right|\, dh  + \frac{24 }{\sqrt{2 \pi^3} T}.
\end{equation}
Set~$T =\varepsilon n^{\gamma}$, where~$\varepsilon>0$, and write
\begin{equation}
  \label{theorem1_proof_rhs_last_ineq}
  \begin{aligned}
    \delta &\le \frac{1}{\pi} \int \limits_{-\varepsilon n^{\gamma}}^{\varepsilon n^{\gamma}} \left|\frac{\varphi_{f,n}(h)- \varphi_{\mathcal{N}}(h)}{h} \right|\, dh \\
    &\quad +\frac{1}{\pi} \int \limits_{-\varepsilon n^{\gamma}}^{\varepsilon n^{\gamma}} \left|\frac{ e^{i h \mexp[n]{S_{f,n}}}  - 1}{h} \right|  \varphi_{\mathcal{N}}(h)\, dh + \frac{24 }{\varepsilon n^{\gamma} \sqrt{2 \pi^3} }\\
    &\le \frac{1}{\pi} \int \limits_{-\varepsilon n^{\gamma}}^{\varepsilon n^{\gamma}} \left|\frac{\varphi_{f,n}(h)- \varphi_{\mathcal{N}}(h)}{h} \right|\, dh + \frac{\sqrt{2}}{\sqrt{\pi}} |\mexp[n]{S_{f,n}}|  + \frac{24 }{\varepsilon n^{\gamma}\sqrt{2 \pi^3} }.
  \end{aligned}
\end{equation}

Finally, we estimate the last integral in~\eqref{theorem1_proof_rhs_last_ineq}. For sufficiently large~$n$ and~$\gamma \le 1/d$, Lemma~\ref{lemma1_growing_arg} yields
\begin{equation}
  \label{th_fin_est_l}
    \int \limits_{-\varepsilon n^{\gamma}}^{\varepsilon n^{\gamma}} \left|\frac{\varphi_{f,n}(h)- \varphi_{\mathcal{N}}(h)}{h} \right|\, dh  \le \frac{C}{n^{1-(d-1)\gamma}} \int \limits_{-\varepsilon n^{\gamma}}^{\varepsilon n^{\gamma}} \varphi_{\mathcal{N}}(h) \, dh \le \frac{\sqrt{2 \pi} C}{n^{1-(d-1)\gamma}},
\end{equation}
where~$C>0$. Notice that due to~\eqref{eq_big_gamma}, this inequality also holds for~$\gamma \in (1/d, 1/(d-1))$.

Collecting all the terms and choosing~$\gamma=1/d$ to attain the best available rate of convergence, we arrive at the desired asymptotic formula
\begin{equation}
  \sup_{n,x}{\left(n^{1/d}|F_{f,n}(x)-F_{\mathcal{N}}(x)|\right)} < +\infty,
\end{equation}
which concludes the proof.

\section*{Acknowledgments}
We are deeply grateful to Christophe Charlier, Tom Claeys, Tamara Grava,  Igor Krasovsky, Arno Kuijlaars, and Oleg Lisovyi for helpful discussions. We also highly appreciate useful remarks and comments of the anonymous reviewers. Our research  is supported by the European Research Council (ERC) under the European Union’s Horizon 2020 research and innovation programme, grant  647133 (ICHAOS). A. B. is also supported by Agence Nationale de Recherche, project ANR-18-CE40-0035, and by the Russian Foundation for  Basic Research, grant 18-31-20031.


\end{document}